\documentclass[12pt]{article}%

\title{A study of Galois objects for algebraic quantum groups}%
\author{Kenny De Commer\footnote{Research Assistant of the Research Foundation - Flanders (FWO -
Vlaanderen).}\\\small Department of Mathematics,  K.U. Leuven\\
\small Celestijnenlaan 200B, B-3001 Leuven, Belgium\\ \\ \small e-mail: kenny.decommer@wis.kuleuven.be}%
\date{}

\usepackage{mathrsfs}%
\usepackage{color}%
\usepackage{amssymb, amsthm, amsfonts, amsxtra, amsmath}%
\usepackage{latexsym}%

\theoremstyle{change}

\newtheorem{Theorem}{Theorem}[section]%
\newtheorem{Def}[Theorem]{Definition}%
\newtheorem{Lem}[Theorem]{Lemma}%
\newtheorem{Prop}[Theorem]{Proposition}%
\newtheorem{Cor}[Theorem]{Corollary}

\begin{document}
\maketitle

\newcommand{\qbino}[2]{\left[\begin{array}{cc}  #1 \\ #2 \end{array}\right]}

\abstract{\noindent We supplement the study of Galois theory for
algebraic quantum groups started in \cite{VDZ1}. We examine the
structure of the Galois objects: algebras equipped with a Galois
coaction such that only the scalars are coinvariants. We show how their
structure is as rich as the one of the quantum groups themselves:
there are two distinguished weak K.M.S. functionals, related by a
modular element, and there is an analogue of the antipode squared. We also show how to reflect the quantum group across
the Galois object to obtain a (possibly) new algebraic quantum
group. We end by considering an example.}

\section{Introduction}

\noindent \textit{Hopf-Galois theory} studies, in a restricted
sense, extensions of unital $k$-algebras $F\subseteq X$ over a field
$k$, with $F$ being the set of coinvariants for a coaction $\alpha:
X\rightarrow X\underset{k}{\otimes} A$ of a Hopf algebra $A$, satisfying the property
that the map $X\underset{F}{\otimes} X\rightarrow
X\underset{k}{\otimes} A: x\underset{F}{\otimes} y \rightarrow
(x\otimes 1)\alpha(y)$ is a bijection. In case $X$ and $A$ are commutative, this can be interpreted geometrically as $X$ being the function space of a bundle
over the spectrum of $F$, with the spectrum of $A$ acting freely and transitively on every fiber.\\


\noindent An interesting situation arises when $F$ reduces to the
scalar field $k$. In this case $X$ is called a \textit{Galois
object}. The most famous instance of this is when $X$ is a field
extension of $k$ and $A$ is the function algebra of a finite group
$G$. The obtained objects are then precisely the finite, normal,
separable field extensions of $k$, with $G$ as Galois group, and the
theory becomes
classical Galois theory.\\



\noindent An important aspect of Galois extensions is that they
provide equivalences of certain categories, see e.g.\!\! \cite{Sch1} for
an overview. As for Galois objects, their isomorphism classes are in
one-to-one correspondence with the isomorphism classes of monoidal
equivalences from the category of corepresentations (i.e. comodules)
of $A$ to a
category of corepresentations of another Hopf algebra $B$.\\


\noindent In the operator algebra framework the same objects appear
under a different name. In \cite{BDV1} a method was constructed to
obtain ergodic coactions of compact quantum groups on C$^*$-algebras, in which irreducible representations
appear with a multiplicity greater than their dimension (which is
impossible for an ordinary compact group). These constructed actions
were termed `of full quantum multiplicity', which is precisely the
condition that the action is Hopf-Galois (when restricting all
C$^*$-algebras to natural dense subalgebras). The
methods however made use of the particular nature of the dual of the
compact quantum groups, which consists of a direct sum of matrix
algebras,
contrasting with the techniques from the Hopf algebra approach. \\

\noindent In this paper we provide a more general way to work
with these objects. Namely, we will consider the structure of Galois
objects for algebraic quantum groups. Algebraic quantum groups were
developed by Van Daele in \cite{VD1}, motivated in turn by finding
the right infinite-dimensional generalization of a
finite-dimensional Hopf algebra which still allows for a dual object
of the same kind, and in providing a purely algebraic framework for
the study of some of the aspects of locally compact quantum groups.
The main differences with a general ordinary Hopf algebra are the
possible lack of a unit in the algebra, and the existence of a
non-zero left invariant functional.\\

\noindent A first study of the general Hopf-Galois theory for
algebraic quantum groups appeared in \cite{VDZ1}, whose main
result is a Morita context between $F$ and the smash product $X\#\hat{A}$ for an
$A$-Galois extension $F\subseteq X$ by an algebraic quantum group
$A$. We refer to that paper for basic notions and results.\\

\noindent We come to the specific content of this paper. \textit{In
the first part}, we study the further algebraic
structure of a Galois object $X$ over an algebraic quantum group
$A$. The main results are the following. There are two distinguished
functionals $\varphi_X$ and $\psi_X$ on $X$, with $\psi_X$ invariant
with respect to the action $\alpha$. They are related by an
invertible element $\delta_X$ inside the multiplier algebra $M(X)$ of $X$, namely $\varphi_X(\,\cdot\,
\delta_X)=\psi_X$. They also satisfy the weak K.M.S. condition: if
$\omega$ is a functional on $X$, we say that it satisfies this
condition if there exists an algebra automorphism $\sigma_\omega$ of
$X$ such that $\omega(y\sigma_\omega(x))=\omega(xy)$ for all $x,y\in
X$. We call a functional satisfying this condition a
\textit{modular} functional. It is interesting to see this
structure, much used in the theory of von Neumann algebras, appear
in a natural way in a purely algebraic setting. There also is a
distinguished automorphism $\theta_X$ on $X$, which plays the
r\^{o}le of the antipode squared. However, there does not seem to be
an analogue of the antipode itself (as is to be expected). Finally, we associate a scaling
constant to $X$, and show that it equals the scaling constant of
$A$. Thus $X$ is almost as rich in structure as $A$. We also show that when working with a $^*$-structure, all these maps become simultaneously diagonizable. We end this
part by considering what happens when $A$ is of a special
type, namely compact or discrete.\\

\noindent \textit{In the second part}, we construct a new algebraic quantum group, starting from a Galois object $X$. We show the connection between Galois objects and monoidal equivalence of module categories, spending some more time on the $^*$-algebraic case. In our calculations, we make essential use of a natural subspace of the dual of $X$. This allows us to make the formulas very transparent.\\

\noindent \textit{In the third part}, we examine a concrete example, which will allow us to construct new examples of algebraic quantum groups of compact type. It illustrates how duality can be used to find in a fairly easy way the structure of a reflected quantum group.\\

\noindent \textit{In the appendix}, we repeat some notions concerning multiplier algebras and algebraic quantum groups.\\

\noindent We now set down conventions and notations. We work over a fixed
field $k$, i.e.\! all algebras are $k$-algebras. Moreover, all
algebras appearing are non-trivial and non-degenerate (which means that
sending an element to the map `left (or right) multiplication with it' is
injective). We denote the multiplier algebra of an algebra $X$ as
$M(X)$, and identify $X$ with its image inside $M(X)$. We also identify $k$ with its image $k\cdot 1$ in $M(X)$. We denote the tensor product over $k$ by $\otimes$. We fix here the
algebraic quantum group $(A,\Delta)$ which we will be working with.
We denote its antipode with $S$ and its counit with $\varepsilon$.
We denote by $\varphi$ a non-zero left invariant functional $\varphi$, positive if
$A$ is a $^*$-algebraic quantum group. As a right invariant
functional $\psi$ we choose $\psi=\varphi\circ S$. We denote the
modular element by $\delta$, and the modular automorphisms of
$\varphi$ and $\psi$ by respectively $\sigma$ and $\sigma'$. We take
the algebraic convention for the coproduct $\hat{\Delta}$ on the dual $\hat{A}$ of $A$: it is determined by $\hat{\Delta}(\omega_1)(x\otimes y) =
\omega(xy)$ for all $x,y\in A$, $\omega_1\in \hat{A}$. The left integral $\hat{\varphi}$ of $\hat{A}$ is given by $\hat{\varphi}(\psi(a\,\cdot\,))=\varepsilon(a)$. For the further
theory of algebraic quantum groups, we refer the reader to
\cite{VD1} and the appendix.



\section{The structure of Galois objects}

\subsection*{\qquad \textit{Definitions}}

\noindent Let $X$ be an algebra. Let $\alpha$ be a \textit{right
coaction} of $A$ on $X$: it is an injective homomorphism
$X\rightarrow M(X\otimes A)$ satisfying
\begin{enumerate} \item [i)]$\alpha(X)(1\otimes A)\subseteq X\otimes A$,
\item [ii)] $(1\otimes A)\alpha(X)\subseteq X\otimes A$, and
\item [iii)] $(\alpha\otimes \iota)\alpha=(\iota\otimes
\Delta)\alpha$,\end{enumerate}

 \noindent where the third property can be made sense of by using the right coverings. The maps
\[\left.\begin{array}{l} T_1: x\otimes a\rightarrow \alpha(x)(1\otimes a),\\ T_2:x\otimes
a\rightarrow (1\otimes a)\alpha(x)\end{array}\right.\] are then bijections, their inverses determined by \[\left.\begin{array}{l} T_1^{-1}: y\otimes S(b)\rightarrow (\iota\otimes S)((1\otimes b)\alpha(y)),\\ T_2^{-1}:y\otimes
S^{-1}(b)\rightarrow (\iota\otimes S^{-1})(\alpha(y)(1\otimes b)).\end{array}\right.\] The injectivity of $\alpha$ implies that $(\iota\otimes \varepsilon)(\alpha(x))=x$ for all $x\in X$.\\


\noindent The coaction is called \textit{reduced} if in addition
$(X\otimes 1)\alpha(X)\subseteq X\otimes A$. The other inclusion
$\alpha(X)(X\otimes1)\subseteq X\otimes A$ follows from this (see the remark after Proposition 2.5. in \cite{VDZ1}).\\

\noindent We can extend $\alpha$ to a map $M(X)\rightarrow M(X\otimes A)$, using the bijections $T_1$ and $T_2$. The \textit{algebra of coinvariants} $F\subseteq M(X)$ is defined as the set of elements $f$ in $M(X)$ such that $\alpha(f)=f\otimes 1$. Remark that it is a unital algebra. The
coaction $\alpha$ is then called \textit{Galois}, if it is reduced and if the map
\[V:X\underset{F}{\otimes} X\rightarrow X\otimes A: x\underset{F}{\otimes} y\rightarrow (x\otimes 1)\alpha(y)\] is
bijective. In fact, the bijectivity already follows from the surjectivity of this map (see Theorem 4.4. in \cite{VDZ1}). Also, the map \[W:X\underset{F}{\otimes} X \rightarrow X\otimes A: x\underset{F}{\otimes} y\rightarrow \alpha(x)(y\otimes 1)\] is
bijective, with the inverse map given by \[W^{-1}(x\otimes a)=
V^{-1}((1\otimes S^{-1}(a))\alpha(x)).\]



\begin{Def} Let $\alpha$ be a right Galois coaction of $A$ on the algebra $X$. Then $(X,\alpha)$ is called a right $A$-\emph{Galois object} if the algebra $F$ of coinvariants reduces to the scalar field $k$.\end{Def}


\noindent In the following $(X,\alpha)$ will be a fixed $A$-Galois
object. For $a\in A$, we denote by $\beta(a)$ the element in $M(X\otimes X)$ which satisfies
\[\left.\begin{array}{l} (x\otimes 1)\beta(a)=V^{-1}(x\otimes a),\\ \beta(a)(1\otimes y) = W^{-1}(y\otimes S(a)).\end{array}\right.\] Using the formula for $W^{-1}$ in terms of $V^{-1}$, it is not difficult to see that $\beta(a)$ is indeed a well-defined multiplier on $X\otimes X$ for each $a\in A$.\\

\noindent We will show later that also the maps
\[x\otimes a\rightarrow \beta(a)(x\otimes 1)\] and \[x\otimes
a \rightarrow (1\otimes x)\beta(a)\] are bijections from $X\otimes
A$ to
$X\otimes X$ (see Corollary \ref{Cor1}).\\


\noindent For
computations we will use the Sweedler notation, denoting $\alpha(x)$ as
$x_{(0)}\otimes x_{(1)}$ (without the summation sign) and $\beta(a)$ as $a^{[1]}\otimes a^{[2]}$. Then by definition we have the identities
\[\left.\begin{array}{l} xa^{[1]}a^{[2]}_{\;\;\;\;(0)}\otimes a^{[2]}_{\;\;\;\;(1)} =
x\otimes a, \\ a^{[1]}_{\;\;\;\;(0)}a^{[2]}x\otimes
a^{[1]}_{\;\;\;\;(1)} = x\otimes S(a), \end{array}\right.\] for all $x\in X, a\in A$. Applying $\iota\otimes \varepsilon$, we obtain the formula \[xa^{[1]}a^{[2]}=\varepsilon(a)x.\] We want to remark and warn however that the use of the Sweedler notation here is more delicate than for Hopf algebras. Indeed, when doing computations with Sweedler notation, it is vital that all expressions are covered. We refer to \cite{DVW1} for a careful analysis of this technique, and to the appendix for the more intuitive approach.\\

\subsection*{\qquad \textit{The existence of the invariant
functionals}}

\noindent For any functional $\omega$ of $X$, we can interpret $(\omega\otimes \iota)\alpha(x)$ in a natural way as a multiplier of $A$. By an invariant functional on $X$ we mean a functional $\omega$ on $X$ such that $(\omega\otimes\iota)\alpha(x)=\omega(x)1$ for all $x\in X$. By a $\delta$-invariant functional we mean a functional $\omega$ on $X$ such that $(\omega\otimes \iota)\alpha(x)=\omega(x)\delta$ for all $x\in X$.

\begin{Prop} There exists a faithful $\delta$-invariant functional $\varphi_X$ on $X$ such that \[(\iota\otimes \varphi)(\alpha(x))=\varphi_X(x)1\] for all $x\in X$.\end{Prop}

\begin{proof} It is clear that the given identity determines $\varphi_X$, since for any $x\in X$, the element $ (\iota\otimes \varphi)(\alpha(x))\in M(X)$ is coinvariant. This map $\varphi_X$
is $\delta$-invariant: for $x,z\in X$ and $a\in A$, we have
\begin{eqnarray*} \varphi_X(x_{(0)})z \otimes x_{(1)}a&=&
x_{(0)}z\otimes\varphi(x_{(1)})  x_{(2)}a\\ &=&x_{(0)}\varphi(x_{(1)})z\otimes
\delta a\\ &=& \varphi_X(x) z\otimes \delta
a.\end{eqnarray*}

\noindent We prove faithfulness. Suppose $x\in X$
is such that $\varphi_X(xy)=0$ for all $y\in X$. Then
\[\varphi(x_{(1)}y_{(1)})x_{(0)}y_{(0)}z=0\qquad \textrm{for all }y,z\in X.\] Using
the Galois property, it follows that \[\varphi(x_{(1)}a)x_{(0)}y=0\qquad \textrm{for all }y\in X \textrm{ and }a\in A. \] The faithfulness of $\varphi$ implies
that $x_{(0)}y\otimes x_{(1)}=0$ for all $y \in X$,
hence $x=0$. Likewise $\varphi_X(yx)=0$ for
all $y\in X$ implies $x=0$.\end{proof}

\begin{Cor} The algebra $X$ has local units:
for any finite set of $x_i\in X$, there exists $y\in X$ with
$yx_i=x_iy=x_i$ for all $i$.\end{Cor}\begin{proof} The proof is the same as the
one
of Proposition 2.6.\! in \cite{DVZ1}, and will be omitted.\end{proof}

\begin{Prop} There exists a non-zero invariant functional $\psi_X$ on $X$.\end{Prop}
\begin{proof}
\noindent Choose $y\in X$ and put
\[\psi_X^y(x)=\varphi_X(x_{(0)}y)\psi(x_{(1)}).\] It is easy to see
that this functional is invariant. Suppose that $\psi_X^y$ is zero
for all $y\in X$. Then $\varphi_X(x)\psi(a)=0$ for all $x\in X$ and $a\in
A$, which is impossible. So we can choose as $\psi_X$ some non-zero
$\psi_X^y$.\end{proof}

\noindent We prove a uniqueness result concerning the invariant functionals. We can follow
the method of Lemma 3.5. and Theorem 3.7. of \cite{VD1} verbatim.

\begin{Prop} If $\psi^{1}_X$ and $\psi^{2}_X$ are two invariant non-zero functionals on $X$, then there exists a scalar $c\in k$ such that $\psi^{1}_X=c\psi^{2}_X$.\end{Prop}
\begin{proof}
First, we show that if $\psi_X$ is a non-zero invariant
functional on $X$, then
\[\{\varphi_X(\,\cdot\, x)\mid x\in X\}=\{\psi_X(\,\cdot\, x)\mid x\in
X\}.\] Choose $x,y,z$ in $X$. Then \[\alpha(xy)(z\otimes 1)=\sum_i
\alpha(xw_i)(1\otimes a_i)\] for some $w_i\in X,a_i\in A$. Likewise
starting with $x,w\in X, a\in A$, there exist $y_i,z_i\in X$ with
\[\sum_i\alpha(xy_i)(z_i\otimes 1)= \alpha(xw)(1\otimes a).\] If we apply
$\psi_X\otimes \varphi$ to these expressions we obtain respectively the equalities
\[\left.\begin{array}{l} \varphi_X(xy)\psi_X(z)=\sum_i \psi_X(xw_i)\varphi(a_i), \\ \sum_i\varphi_X(xy_i)\psi_X(z_i)= \psi_X(xw)\varphi(a). \end{array}\right.\] Choosing either $z$ with $\psi_X(z)=1$ or
$a$ with $\varphi(a)=1$, we get respectively $\subseteq$ and $\supseteq$.\\

\noindent Suppose now that $\psi_X^1$ and $\psi_X^2$ are invariant. Choose
$y,z_1\in X$ with $\varphi_X(yz_1)=1$ and take $z_2\in X$ with $\psi_X^1(\cdot
z_1)=\psi_X^2(\cdot z_2)$. Choosing $x\in X$, applying $\psi_X^i\otimes \varphi$ to
$(x\otimes 1)\alpha(yz_i)$ and writing this last expression as
$\sum_j (1\otimes a_j)\alpha(w_jz_i)$ for certain $w_j\in X, a_j\in A$, we see that
$\psi_X^1(x)=\varphi_X(yz_2)\psi_X^2(x)$, proving that all invariant
functionals are scalar multiples of each other.\end{proof}

\subsection*{\qquad \textit{The existence of the modular element}}

\noindent Let $\psi_X$ be a non-zero invariant functional on $X$. We prove the existence of a modular element $\delta_X$,
relating the functionals $\varphi_X$ and $\psi_X$. We first prove a lemma:

\begin{Lem}\label{lem3} For all $x,y\in X$ we have \[S((\psi_X\otimes \iota)((x\otimes 1)\alpha(y)))=(\psi_X\otimes \iota)(\alpha(x)(y\otimes 1)).\] \end{Lem} \begin{proof}
\noindent Choose $a\in A$ and $x,y\in X$. Pick $z_i \in X$ and $b_i\in A$ such that \[(1\otimes a)\alpha(y)= \sum_i z_i\otimes b_i.\] Then by the formula for $T_2^{-1}$ given above, we have \[y\otimes S(a)=\sum_i\alpha(z_i)(1\otimes S(b_i)).\] If we denote $w=S((\psi_X\otimes \iota)((x\otimes 1)\alpha(y)))$, then
 \begin{eqnarray*} w S(a) &=&  \sum_i \psi_X(xz_i)S(b_i)\\ &=& \sum_i \psi_X(x_{(0)}z_{i(0)})x_{(1)}z_{i(1)}S(b_i)\\ &=&\psi_X(x_{(0)}y)x_{(1)}S(a).\end{eqnarray*} Since $a$ was arbitrary, the formula is proven.

\end{proof}

\begin{Prop} There exists a unique invertible element $\delta_X\in M(X)$ such that $\varphi_X(x\delta_X)=\psi_X(x)$
for all $x\in X$.\end{Prop}

\begin{proof} We show first that for all $x\in X$:
\[\psi_X(x)=0 \qquad\Rightarrow\qquad \psi(x_{(1)})x_{(0)}=0.\] We
know that $\psi_X(x)=0$ implies $\psi_X^w(x)=0$
for all $w\in X$, i.e. $\psi(x_{(1)})\varphi_X(x_{(0)}w)=0$ for all $w\in X$.
So $\psi(x_{(1)})x_{(0)}=0$ by the faithfulness of $\varphi_X$.\\

\noindent This means that $\psi(x_{(1)})x_{(0)}=c_x\delta_X'$ for some
multiplier $\delta_X'\in M(X)$ and some number $c_x\in k$. Now
$x\rightarrow c_x$ is easily seen to be a non-zero invariant functional, and replacing $\psi_X$ by this invariant functional, we obtain
$\psi(x_{(1)})x_{(0)}=\psi_X(x)\delta_X'$.\\

\noindent Now we show that $\delta_X'$ has an inverse $\delta_X$, and that
$\varphi_X(x\delta_X)=\psi_X(x)$. Choose $y\in X$ with $\psi_X(y)=1$,
then
\begin{eqnarray*}
\psi_X(x\delta_X')&=&\psi_X(xy_{(0)})\varphi(S(y_{(1)}))\\ &=&
\psi_X(x_{(0)}y)\varphi(x_{(1)})\\&=&\varphi_X(x).\end{eqnarray*}
Since furthermore $\{\varphi_X(\,\cdot\, x)\}=\{\psi_X(\,\cdot\, x)\}$, we
have that for any $x\in X$ there exists $y\in X$ with $y\delta_X'=x$. To show
that also left multiplication is surjective, we use another
argument. Take $x\in X$ and $a\in A$ with $\psi(a)=1$. Write
$x\otimes a$ as $\sum_i p_{i(0)}q_i\otimes p_{i(1)}$ for certain $p_i,q_i\in X$, and put $y=\sum_i \psi_X(p_i)q_i$.
Then $\delta_X'y= \sum_i \psi(p_{i(1)})p_{i(0)}q_i=\psi(a)x=x$.\\

\noindent Hence we obtain the formula
\[\psi_X(x)=\varphi_X(x\delta_X) \qquad \textrm{for all } x\in X.\] By the faithfulness of $\varphi_X$, this uniquely determines $\delta_X$.\end{proof}

\begin{Cor} Any invariant non-zero functional is faithful.\end{Cor}

\subsection*{\qquad \textit{The modularity of the invariant
functionals}}

\noindent We first prove some identities.

\begin{Lem}\label{lem1}   For all
$x\in X$ and $a\in A$, we have
\begin{enumerate}\item [i\textrm{)}]$
\varphi(ax_{(1)})x_{(0)}=\varphi_X(a^{[2]}x)a^{[1]}$,
\item [ii\textrm{)}]$
\varphi(x_{(1)}S(a))x_{(0)}=\varphi_X(xa^{[1]})a^{[2]}$.\end{enumerate}\end{Lem}\begin{proof}
The first equation follows from \begin{eqnarray*}
\varphi(ax_{(1)})zx_{(0)}&=&
\varphi(a^{[2]}_{\;\;\;(1)}x_{(1)})za^{[1]}a^{[2]}_{\;\;\;(0)}x_{(0)}\\
&=& \varphi_X(a^{[2]}x)za^{[1]},\end{eqnarray*} for all $a\in A$ and
$x,z\in X$. The second follows from \begin{eqnarray*}
\varphi(x_{(1)}S(a))x_{(0)}z&=&
\varphi(x_{(1)}a^{[1]}_{\;\;\;(1)})x_{(0)}a^{[1]}_{\;\;\;(0)}a^{[2]}z\\
&=& \varphi_X(xa^{[1]})a^{[2]}z,\end{eqnarray*} for all $a\in A$ and
$x,z\in X$.\end{proof}

\begin{Lem} For all $y,p,q\in X$ and $a\in A$, we have \[\varphi_X(a^{[2]}y)\varphi_X(pa^{[1]}q)=\varphi_X(yb^{[1]})\varphi_X(pb^{[2]}q),\qquad\textrm{where } b=(S^{-1}\sigma)(a).\]\end{Lem}
\begin{proof}

 Using the identities of the previous lemma, we get \begin{eqnarray*} \varphi_X(yb^{[1]})\varphi_X(pb^{[2]}q)&=& \varphi(y_{(1)} \sigma(a))\varphi_X(py_{(0)}q)\\ &=& \varphi(ay_{(1)})\varphi_X(py_{(0)}q)\\ &=& \varphi_X(a^{[2]}y)\varphi_X(pa^{[1]}q).\end{eqnarray*}
\end{proof}

\noindent We show now that $\varphi_X$ is modular.

\begin{Prop}
There exists an automorphism $\sigma_X$ of $X$ such that
\[\varphi_X(y\sigma_X(x))=\varphi_X(xy)\qquad \textrm{for all
}x,y\in X.\] Furthermore, $\varphi_X\circ \sigma_X=\varphi_X$.\end{Prop}\begin{proof} Choose $x\in X$, and write $x$ as a sum of elements of the form $\varphi_X(pa^{[1]}q)a^{[2]}$ with $p,q\in X$ and $a\in A$. Define $w$ as $\sum \varphi_X(pb^{[2]}q)b^{[1]}$ with $b=(S^{-1}\sigma)(a)$. Then the previous lemma shows that $\varphi_X(yw)=\varphi_X(xy)$ for all $y\in X$.\\


\noindent It is clear that $w$ is uniquely determined by this
property, so we can denote $w=\sigma_X(x)$. Standard arguments imply
that $\sigma_X$ is indeed an algebra automorphism. It will leave
$\varphi_X$ invariant because $X^2=X$. \end{proof}

\begin{Cor} The functional $\psi_X$ is modular with modular
automorphism \[\sigma_X'(x)=\delta_X\sigma_X(x)\delta_X^{-1}.\]\end{Cor}

\subsection*{\qquad \textit{Formulas}}

\noindent In this section and the next, we collect some formulas. They strongly resemble the formulas which hold in algebraic quantum groups, and also their proofs are mostly straightforward adaptations. Nevertheless, it is remarkable that Galois objects carry as rich a structure as the algebraic quantum groups themselves.

\begin{Prop}\label{prop1} For all $a\in A$, we have
\begin{enumerate} \item [i)] $\alpha\circ \sigma_X'=(\sigma_X'\otimes
S^{-2})\circ \alpha$,
\item [ii)] $((S^{-1}\sigma)(a))^{[1]}\otimes
((S^{-1}\sigma)(a))^{[2]}=\sigma_X(a^{[2]})\otimes a^{[1]}$.
\end{enumerate}

\end{Prop}

\begin{proof} Choose $x,y\in X$ and $a\in A$. Then \begin{eqnarray*} (\psi_X\otimes \varphi)((y\otimes a)\alpha(\sigma_X'(x)))&=& \psi_X(y_{(0)}\sigma_X'(x))\varphi(aS^{-1}(y_{(1)}))\\ &=& \psi_X(xy_{(0)}) \varphi(aS^{-1}(y_{(1)}))\\ &=& \psi_X(x_{(0)}y)\varphi(aS^{-2}(x_{(1)}))\\ &=& \psi_X(y\sigma'(x_{(0)}))\varphi(aS^{-2}(x_{(1)})),\end{eqnarray*}
applying Lemma \ref{lem3} twice. As $\varphi$ and $\psi_X$ are faithful, the first identity follows. The second formula was essentially proven in the previous section.\end{proof}


\begin{Cor}\label{Cor1}  The maps
\[\left.\begin{array}{l} x\otimes a\rightarrow \beta(a)(x\otimes 1),\\x\otimes
a\rightarrow (1\otimes x)\beta(a)\end{array}\right.\] are bijections from $X\otimes
A$ to $X\otimes X$\end{Cor}\begin{proof} \noindent This follows from the second formula.\end{proof}  \noindent Note that this fact is not clear at first sight. It allows us for example to construct the
\textit{Miyashita-Ulbrich} action in this context: $A$ acting on the
right of $X$ by $x\cdot a= a^{[1]}xa^{[2]}$, making it a Yetter-Drinfel'd module. It also allows us to regard $\beta$ rather as a map $\tilde{\beta}:A\rightarrow M(X^{\textrm{op}}\otimes X)$. Indeed: if we are ignorant of the previous corollary, we would only know that $\tilde{\beta(a)}$ is a left multiplier of $X^{\textrm{op}}\otimes X$. Now as is the case for Galois objects over Hopf algebras, the map $\tilde{\beta}$ will be a homomorphism. The argument for this is simple: choose $x\in X$ and $a,b\in A$, and write $xb^{[1]}\otimes b^{[2]}= \sum_i p_i \otimes q_i$ for certain $p_i,q_i\in X$. Then $\sum_i (p_i\otimes a)\alpha(q_i) = x\otimes ab$. Applying $V^{-1}$, we obtain $ \sum_i p_ia^{[1]}\otimes a^{[2]}q_i = x(ab)^{[1]}\otimes (ab)^{[2]}$, so $xb^{[1]}a^{[1]}\otimes a^{[2]}b^{[2]} = x(ab)^{[1]}\otimes (ab)^{[2]}$. This proves that $\tilde{\beta}$ is a homomorphism.\\

\noindent The following proposition collects some formulas concerning the modular elements.
\begin{Prop} The following identities hold:
\begin{enumerate}
\item [iii)] $\alpha(\delta_X)=\delta_X\otimes \delta$,
\item [iv)] $\tilde{\beta}(\delta) = \delta_X^{-1}\otimes \delta_X$,
\item [v)] $\sigma_X(\delta_X)=\tau^{-1}\delta_X$,
\end{enumerate} where in the last formula $\tau$ denotes the scaling
constant of $A$. \end{Prop}
\begin{proof} First note that for any $x,y\in X$ we have $\varphi_X(xy_{(0)}y_{(1)} = \varphi_X(x_{(0)}y)S^{-1}(x_{(1)})\delta$, which is proven in the same way as \ref{lem3}, using the $\delta$-invariance of $\varphi_X$. Then if $x,y\in X$, we have \begin{eqnarray*} \varphi_X(x(y\delta_X)_{(0)})(y\delta_X)_{(1)} &=& \varphi_X(x_{(0)}y\delta_X)S^{-1}(x_{(1)})\delta \\ &=& \psi_X(x_{(0)}y)S^{-1}(x_{(1)})\delta \\ &=& \psi_X(xy_{(0)})y_{(1)}\delta \\ &=& \varphi_X(xy_{(0)}\delta_X)y_{(1)}\delta.\end{eqnarray*} By faithfulness of $\varphi_X$ we have $\alpha(y\delta_X)=\alpha(y)(\delta_X\otimes \delta)$, hence $\alpha(\delta_X)=\delta_X\otimes \delta$ by definition of $\alpha$ on $M(X)$.\\

\noindent For the second formula, we have to prove that $x(a\delta)^{[1]} \otimes (a\delta)^{[2]}= x\delta_X^{-1}a^{[1]} \otimes a^{[2]}\delta_X $ for all $a\in A$ and $x\in X$. This follows immediately by applying $V$ and using the previous formula.\\

\noindent As for the final formula, we have for any $x\in X$ that \begin{eqnarray*} \varphi_X(\delta_Xx) &=& \varphi(\delta x_{(1)})\delta_X x_{(0)}\\ &=& \tau^{-1} \varphi(x_{(1)}\delta)\delta_X (x_{(0)}\delta_X)\delta_X^{-1}\\ &=& \tau^{-1}\varphi((x\delta_X)_{(1)})\delta_X(x\delta_X)_{(0)}\delta_X^{-1}\\ &=& \tau^{-1}\varphi_X(x\delta_X),\end{eqnarray*} which means exactly that $\sigma_X(\delta_X)=\tau^{-1}\delta_X$.\end{proof}

\begin{Cor} If $\varphi_X'$ is another $\delta$-invariant functional, then there exists $c\in k$ with $\varphi_X'=c\varphi_X$.\end{Cor}
\begin{proof}
This follows immediately by the fact that $\varphi_X'(\,\cdot\, \delta_X)$ is invariant.\end{proof}

\subsection*{\qquad \textit{The square of no antipode}}

\noindent On $X$ there is a natural unital left $\hat{A}$-module algebra
structure defined by $\omega_1\cdot x= (\iota\otimes
\omega_1)\alpha(x)$ for $x\in X$ and $\omega_1\in \hat{A}$. The unitality means that $\hat{A}\cdot X=X$. It allows us to extend the action to a left action of $M(\hat{A})$ on $X$.\\

\noindent Consider then the map \[\theta_X: X\rightarrow X:
x\rightarrow \sigma_X(\hat{\delta}\cdot x),\] where $\hat{\delta}$ is the modular element of $\hat{A}$. It is
a bijective homomorphism, which can be shown using the module algebra structure and the fact that
$\hat{\delta}$ is grouplike. This $\theta_X$ plays the r\^{o}le of `the square of the
antipode' for $X$, even though there is no `antipode of $X$'. Indeed: in case $X=A$ and $\alpha= \Delta$, then $\theta_X$ is exactly $S^2$.
We can use $\theta_X$ to complete our set of formulas.

\begin{Prop}
For all $x\in X$:
\begin{enumerate}\item[vi)] $\alpha\circ \sigma_X = (\theta_X\otimes
\sigma)\circ \alpha$, \item[vii)]$\alpha\circ \theta_X =
(\theta_X\otimes S^2)\circ \alpha$ , \item[viii)] $\alpha\circ
\theta_X = (\sigma_X\otimes \sigma'^{-1})\circ \alpha$,
\item [ix)] $\sigma_X\circ \theta_X=\theta_X\circ \sigma_X$,
\item [x)] $\theta_X(\delta_X)=\delta_X$,
\item[xi)] $\varphi_X\circ \theta_X = \varphi_X(\delta_X^{-1}\cdot \delta_X)= \tau \varphi_X$.
\end{enumerate} \end{Prop} \begin{proof} We first make some remarks. The vector space $M(\hat{A})$ can still be seen as a subspace of the vector space dual $A^*$ of $A$. As such the element $\hat{\delta}$ corresponds with $\varepsilon \circ \sigma^{-1}=\varepsilon \circ \sigma'^{-1}$ (see Proposition 5.14. in \cite{Kus2}). Thus we can also write $\theta_X(x)= \varepsilon(\sigma^{-1}(x_{(1)})) \sigma_X(x_{(0)})$. Also, the above formulas are known to hold in case $X=A$ and $\alpha=\Delta$. We will use them in the course of the proof.\\

\noindent Take $x,y\in X$ and $a\in A$. Then \begin{eqnarray*} \varphi_X(y\theta_X(x_{(0)}))\varphi(a\sigma(x_{(1)}))&=& \varphi_X((\hat{\delta}\cdot x_{(0)})y)\varphi(x_{(1)}a)\\ &=& \varphi_X(x_{(0)}y)\varepsilon (\sigma'^{-1}(x_{(1)})) \varphi(x_{(2)}a)\\ &=&\varphi_X(x_{(0)}y)\varphi(\sigma'^{-1}(S^{-2}(x_{(1)}))a)\\ &=& \varphi_X(xy_{(0)})\psi(aS^{-1}(y_{(1)}))\\ &=& \varphi_X(y_{(0)}\sigma_X(x))\psi(aS^{-1}(y_{(1)}))\\ &=& \varphi_X(y\sigma_X(x)_{(0)})\varphi(a\sigma_X(x)_{(1)}),\end{eqnarray*} which proves the equality in $vi)$. The equality in $vii)$ then follows by the previous one, and the fact that $\theta_A=S^2$ in case $X=A$ and $\alpha=\Delta$.\\

\noindent Further, \begin{eqnarray*} \varphi_X(y\theta_X(x_{(0)}))\psi(S^2(x_{(1)})a) &=& \varphi_X(y\sigma_X(x_{(0)}))\hat{\delta}(x_{(1)})\psi(S^2(x_{(2)})a)\\ &=& \varphi_X(y\sigma_X(x_{(0)}))\varepsilon (\sigma'^{-1}(x_{(1)}))\psi(S^{2}(x_{(2)})a)\\ &=& \varphi_X(y\sigma_X(x_{(0)}))\psi(\sigma'^{-1}(x_{(1)})a),\end{eqnarray*} which together with $vi)$ proves $viii)$. The commutation in $ix)$ is clear. \\

\noindent As for $x)$ we have $\theta_X(\delta_X) = \sigma_X(\delta_X)\varepsilon(\sigma^{-1}(\delta))$, which equals $\delta_X$ by the formula $v)$. The same formula $v)$ also shows immediately the validity of $xi)$. This concludes the proof.
\end{proof}

\noindent We separate formulas which give a more direct connection
with the defining property of an antipode.

\begin{Prop}\label{prop2} For all $x,y\in X$ and $a\in A$, we have \begin{enumerate}\item[xiii)] $S(a)^{[1]}\otimes S(a)^{[2]}= \theta_X(a^{[2]})\otimes a^{[1]}$,\item[xiv)] $x_{(1)}^{\;\;\;[1]}\theta_X(x_{(0)})y\otimes x_{(1)}^{\;\;\;[2]} = y\otimes x$,\item[xv)] $\theta_X(a^{[2]})a^{[1]}y=\varepsilon(a)y$, \item[xvi)] $S^2(a)^{[1]}\otimes S^2(a)^{[2]}= \theta_X(a^{[1]})\otimes \theta_X(a^{[2]})$. \end{enumerate}\end{Prop}
\begin{proof}

\noindent Applying $(\iota\otimes \varphi_X(\,\cdot\, x))$ to
$\theta_X(a^{[2]})\otimes a^{[1]}$ and using formula $vi)$, we get
\begin{eqnarray*} \varphi_X(a^{[1]}x)\theta_X(a^{[2]})&=&
\varphi_X(\sigma_X^{-1}(x)a^{[1]})\theta_X(a^{[2]})\\ &=&
\varphi(\sigma_X^{-1}(x)_{(1)}S(a))\theta_X(\sigma_X^{-1}(x)_{(0)})\\&=&
\varphi(\sigma^{-1}(x_{(1)})S(a))x_{(0)}\\ &=&
\varphi(S(a)x_{(1)})x_{(0)}\\ &=& \varphi_X(S(a)^{[2]}x)
S(a)^{[1]},\end{eqnarray*} which gives the first formula.\\

\noindent As for the second formula, we have to show that for
all $z\in X$ we have \[\varphi_X(x_{(1)}^{\;\;\;[2]}z)
x_{(1)}^{\;\;\;[1]}\theta_X(x_{(0)})y=\varphi_X(xz)y.\] This reduces to proving that \[\varphi(x_{(1)}z_{(1)})y_{(0)}\theta_X(x_{(0)})y=\varphi_X(xz)y.\] This follows again by formula $vi)$ and the defining property of $\varphi_X$. \\

\noindent The last formulae are a direct consequence of the
first.\end{proof}


\subsection*{\qquad \textit{Concerning $^*$-structures}}

\noindent We now look at the case $k=\mathbb{C}$, $A$ a
$^*$-algebraic quantum group and $X$ a $^*$-algebra such that
$\alpha$ is $^*$-preserving. It is easy to see then that $\varphi_X$ is a positive functional, by which we mean that $\varphi_X(x^*x)\geq 0$ for all $x\in X$. We show that $\psi_X$ is positive. As for $^*$-algebraic quantum groups, this is a non-trivial statement. In that case, the first proof of this statement consisted of establishing an analytic structure on the $^*$-algebraic quantum group (see \cite{Kus1}). In \cite{DCV1}, we found an easier way to arrive at this, in the meantime showing something more about the structure of $^*$-algebraic quantum groups. Namely, almost all structure maps on $A$ are
diagonizable, i.e. there exists a basis of simultaneous
eigenvectors for $\sigma$, $S^2$ and left and right
multiplication with $\delta$, which implies for example that the scaling constant is 1. We prove now that also all structure maps on $X$ are diagonizable.\\

\noindent For instance, take $x\in X$ and choose $w\in X$ with $\varphi_X(w)=1$. Write $x\otimes w$ as a sum of $ya^{[1]}\otimes a^{[2]}$ for certain $y\in X$ and $a\in A$. Write $a=\sum a_i$ with the $a_i$ eigenvectors for left multiplication with $\delta$. Then \begin{eqnarray*} x\delta_X^n &=& \sum \varphi_X(a^{[2]}) ya^{[1]}\delta_X^n \\ &=& \sum\varphi_X(\delta_X^{n}(\delta^{-n}a)^{[2]})y(\delta^{-n}a)^{[1]}\\ &\in & \textrm{Span}\{\omega(a_i^{[2]})ya_i^{[1]}\mid \omega\in X^*\},\end{eqnarray*} showing that $\textrm{Span}\{x\delta_X^n\mid n\in \mathbb{Z}\}$ is finite-dimensional. The same technique shows that $\textrm{Span}\{\delta_X^nx\mid n\in \mathbb{Z}\}$ and $\textrm{Span}\{\theta_X^n(x)\mid n\in \mathbb{Z}\}$ are finite-dimensional.  Since the left action of $\hat{\delta}$ is diagonizable, we also have that  $\textrm{Span}\{\sigma_X^n(x)\mid n\in \mathbb{Z}\}$ is finite-dimensional. As all these operations are easily shown to be self-adjoint with respect to the scalar product $\langle x,y\rangle =\varphi_X(y^*x)$ on $X$, and since they commute, we obtain:

\begin{Theorem} There exists a basis $\{x_i\}$ of $X$ such that all $x_i$ are eigenvectors for $\sigma_X, \theta_X$ and multiplication with $\delta_X$ to the left and to the right.\end{Theorem}

\noindent It follows then also that $\delta_X$ is in
fact of the from $(\delta_X^{1/2})^2$ for some self-adjoint element
$\delta_X^{1/2}\in M(X)$. In particular, choosing $z\in X$ with
$\varphi_X ((\delta_X^{-1/2}z)^*\delta_X^{-1/2}z)=1$, we have for
any $x\in X$ that
\begin{eqnarray*} \psi_X(x^*x) &=&
\psi_X(x^*x)\varphi_X(z^*\delta_X^{-1}z)\\&=&
\varphi_X(z^*x_{(0)}^*x_{(0)'}z)\psi(x_{(1)}^*x_{(1)'})\\ &\geq &
0,\end{eqnarray*}  showing

\begin{Cor} The functional $\psi_X$ is positive.\end{Cor}

\noindent We also have to say something about how $\beta$ and $^*$ correspond.
\begin{Prop} For all $a\in A$, we have \[\beta(a)^*= (\Sigma\beta)(S(a)^*),\] where $\Sigma$ denotes the flip map.\end{Prop}\begin{proof} For any $x\in X$, $a\in A$, we have \[\varphi(ax_{(1)})x_{(0)}=\varphi_X(a^{[2]}x)a^{[1]}.\] Applying $^*$, we see that \[\varphi(x_{(1)}^*S(S(a)^*))x_{(0)}^*=\varphi_X(x^*a^{[2]*})a^{[1]*}.\] Since the left hand side equals $\varphi_X(x^*(S(a)^*)^{[1]})(S(a)^*)^{[2]}$, the formula is proved.\end{proof}


\subsection*{\qquad \textit{Concerning compactness and
discreteness}}

\begin{Def}\label{def1} A nondegenerate algebra $X$ is called \emph{of compact type} if $X$ has a unit. It is called \emph{of discrete type} if every subspace of the form $xX$ or $Xx$, with
$x\in X$, is finite dimensional.\end{Def}

\noindent \textit{Remark:} This terminology is not standard, and we use it solely in this subsection.\\

\begin{Theorem} The algebra $X$ is of compact type iff $A$ is an algebraic quantum group of compact type. The algebra $X$ is of discrete type iff $A$ is an algebraic quantum group of discrete type.
\end{Theorem}


\begin{proof} If $A$ is compact,
then $\alpha(x)\in X\otimes A$ for any $x\in X$. Choosing $x\in X$
with $\varphi_X(x)=1$, we have that $(\iota\otimes
\varphi)\alpha(x)\in X$ is a unit of $X$.\\

\noindent If $X$ is compact, then $(x\otimes 1)\alpha(1)\in X\otimes A$ for
all $x\in A$. Choosing $x$ with $\varphi_X(x)=1$, it is again easy
to see that $(\varphi_X\otimes \iota)((x\otimes 1)\alpha(1))\in A$
is a
unit in $A$.\\

\noindent Now suppose that $A$ is an algebraic quantum group of discrete type. Choose a non-zero left
cointegral $h\in A$, so $ah=\varepsilon(a)h$ for all $a\in A$. We
can scale $h$ so that $\varphi(h)=1$. Then for all $x,y\in X$, we
have $\varphi_X(xh^{[1]})h^{[2]}y = \varphi(x_{(1)}h)x_{(0)}y= xy$.
This shows that $Xy$ is finite dimensional. Also $yX$ is finite
dimensional, by a similar reasoning.\\

\noindent Conversely, suppose that $X$ is an algebra of discrete type. Take $a\in A$ and $x\neq
0$ fixed in $X$. Write $xa^{[1]}\otimes a^{[2]}$ as $\sum_i
p_i\otimes q_i$, and choose $y\in X$ such that $p_iy=p_i$ for all
$i$. Then
\begin{eqnarray*} \textrm{dim}\{Aa\} &=& \textrm{dim}\{x\otimes ba
\mid b\in A\}\\ &=& \textrm{dim}\{\sum_i p_i yb^{[1]}\otimes
b^{[2]}q_i\mid b\in A\} \\ &\leq & \textrm{dim span}\{\sum_i p_i
w\otimes zq_i \mid w,z\in X\}\\ &<& \infty.\end{eqnarray*} We show
that this is sufficient to conclude that $A$ is an algebraic quantum group of discrete type.\\

\noindent First, applying $S$, we see that also all $aA$ are finite
dimensional. Choose $a\in A$ with $\varepsilon(a)=1$. Write $I=AaA$,
which is a finite-dimensional ideal. Because $\varphi$ is faithful, we can choose some
$\omega=\varphi(\cdot b)\in \hat{A}$ such that
$\omega_{|I}=\varepsilon_{|I}$. Take $e\in A$ with $ae = a$. Then
for all $x\in A$, we have
\begin{eqnarray*} \varphi(xab) &=& \varphi(xaeb)\\ &=& \omega_{|I}(xae)
\\ &=& \varepsilon(x).\end{eqnarray*} Hence
$\varepsilon\in \hat{A}$, and
$A$ is an algebraic quantum group of discrete type.\end{proof}

\noindent Note that the proof above shows that the terminology we used is consistent: an algebraic quantum group is of discrete type (in the sense of \cite{VD1}) iff its underlying algebra is of discrete type (as defined in Definition \ref{def1}). Also note that if $k=\mathbb{C}$ and $X$ is a $^*$-algebra, the
condition `$X$ is of discrete type' is equivalent with $X$ being a direct
sum of finite-dimensional matrix algebras.



\section{Reflecting an algebraic quantum group across a Galois
object}


\noindent It is known that if $A$ is of compact type (or more generally a
Hopf algebra) and $X$ is unital, then a second Hopf algebra $C$ can
be constructed from $A$ and $X$. Moreover, this $C$ has a left
coaction on $X$, and $X$ will be what is termed a $C$-$A$-bi-Galois
object: it is at the same time a left $C$-Galois object and right
$A$-Galois object, with the two coactions commuting. We show that
the same holds in our setting. Our approach is based on duality: we first construct $\hat{C}$, which is in fact more natural.\\

\noindent We make a remark concerning notation. We will always denote elements in the dual $X^*$ of $X$ by $\omega,\omega',\omega'',\ldots$ and the elements in the dual $A^*$ of $A$ by $\omega_1,\omega_2,\omega_3,\ldots$ When elements of $X^*$ are indexed by some set $I$, we will put the index in superscript, so then these elements take the form $\omega^i$.\


\subsection*{\qquad \textit{The dual of $X$}}

\noindent \begin{Def} The \emph{restricted dual} of $X$ is the vector space $\hat{X}=\{\varphi_X(\,\cdot\, x)\mid x\in X\}$ inside the dual
$X^*$ of $X$.\end{Def}

\noindent We have shown that
\begin{eqnarray*}\hat{X}&=&\{\varphi_X(x\,\cdot\, )\mid x\in X\}\\&=&
\{\psi_X(\,\cdot \,x)\mid x\in X\}\\ &=& \{\psi_X(x\,\cdot\, )\mid x\in
X\}.\end{eqnarray*}

\noindent  The left $\hat{A}$-module structure on $X$
leads to a \textit{right} $\hat{A}$-module structure on $X^*$, which
restricts to a natural right $\hat{A}$-module structure on
$\hat{X}$. More concretely, we have
\[\psi_X(\,\cdot\, x)\cdot \omega_1 = \omega_1(S(x_{(1)}))\psi_X(\,\cdot\, x_{(0)}),\qquad
\textrm{for all } x\in X,\omega_1\in \hat{A},\] by using Lemma \ref{lem3}.\\

\noindent We can dualize the multiplication on $X$ to a map
\[\Delta_X:\hat{X}\rightarrow (X\otimes X)^*.\] We denote the image of
$\omega$ by $\omega^{(1)}\otimes \omega^{(2)}$. While we can not say
that this element is `in $M(\hat{X}\otimes \hat{X})$', since $\hat{X}$
has no multiplication, we do have that expressions such as
$\omega^{(1)}\otimes (\omega^{(2)}\cdot \omega_1)$ with $\omega\in
\hat{X}$ and $\omega_1\in \hat{A}$, are again elements of
$\hat{X}\otimes \hat{X}$, and that this provides bijections between
$\hat{X}\otimes \hat{A}$ and $\hat{X}\otimes \hat{X}$. For example,
the map \[V^t:\hat{X}\otimes \hat{A}\rightarrow \hat{X}\otimes
\hat{X}: \omega\otimes \omega_1\rightarrow \omega^{(1)}\otimes
(\omega^{(2)}\cdot \omega_1)\] is just the dual of the map $V$: \[
(V^t (\omega\otimes \omega_1))(x\otimes y)=(\omega\otimes
\omega_1)(V(x\otimes y)).\] Again in more concrete terms, we have $V^t((\varphi_X(\,\cdot\, x)\otimes \varphi(a\,\cdot\,))= \varphi_X(\,\cdot\, a^{[1]}x)\otimes \varphi_X(a^{[2]}\,\cdot\,)$ for $x\in X,a\in A$, by using Lemma \ref{lem1}.\! $i)$.\\

\noindent The space $\hat{X}$ also carries a \textit{natural} $\hat{A}$-valued $k$-bilinear form, determined by \[\lbrack\,\omega,\omega'\,\rbrack_{\hat{A}}(a)=(\omega\otimes\omega')(\beta(a)).\] If we then let $\hat{A}$ act on the left of $\hat{X}$ by $\omega_1\cdot \omega := \omega\cdot \hat{S}^{-1}(\omega_1)$ for $\omega\in \hat{X}$ and $\omega_1\in \hat{A}$, we have \begin{Prop} The form $\lbrack\,\cdot\,,\,\cdot\,\rbrack_{\hat{A}}$ is $\hat{A}$-bilinear, i.e. for all $\omega,\omega'\in \hat{X}$ and $\omega_1\in \hat{A}$, \[\lbrack \,\omega,\omega'\cdot \omega_1\,\rbrack_{\hat{A}}=\lbrack \,\omega,\omega'\,\rbrack_{\hat{A}}\cdot \omega_1,\] \[\lbrack\, \omega_1\cdot\omega,\omega' \,\rbrack_{\hat{A}}=\omega_1\cdot\lbrack \,\omega,\omega'\,\rbrack_{\hat{A}} .\]
\end{Prop} \begin{proof} This is just a reformulation of the formulas 2.1.2. and 2.1.3. of Lemma 2.7. in \cite{Sch1}. To proof the linearity on the right for example, note that if $\omega= \varphi_X(x\,\cdot\,)$, then \[\lbrack\,\omega,\omega'\,\rbrack_{\hat{A}}(a)=\varphi(x_{(1)}S(a))\omega'(x_{(0)}),\] using the formulas of Lemma 2.8. Then \begin{eqnarray*} (\lbrack\, \omega,\omega'\,\rbrack_{\hat{A}}\cdot \omega_1)(a) &=& \varphi(x_{(1)}S(a_{(1)}))\omega'(x_{(0)})\omega_1(a_{(2)})\\ &=& \varphi (x_{(2)}S(a))\omega'(x_{(0)})\omega_1(x_{(1)})\\ &=& \lbrack \,\omega,\omega'\cdot \omega_1\,\rbrack_{\hat{A}}(a).\end{eqnarray*} The other identity can be proven in the same way.\end{proof} \noindent The following formula shows how the bracket behaves with respect to the left action: \begin{Lem}\label{lem5} For $\omega,\omega',\omega''\in \hat{X}$, we have \[\lbrack \,\omega,\omega'\,\rbrack_{\hat{A}}\cdot \omega''= \omega''\cdot \lbrack \,\theta_X^{-1}(\omega'),\omega\,\rbrack_{\hat{A}},\] with $\theta_X^{-1}(\omega')=\omega'\circ \theta_X^{-1}$. \end{Lem}\begin{proof} This follows from formula $xiii)$ of Proposition \ref{prop2}.
\end{proof}

\subsection*{\qquad \textit{Construction of the reflected algebraic quantum group}}

\noindent On $\hat{X}$ we can consider the vector space $B$ spanned by the operators $\lbrack \,\omega,\omega'\,\rbrack_B$ defined by the identity \[\lbrack\,  \omega,\omega'\,\rbrack_B \cdot \omega''  = \omega\cdot \lbrack \,\omega',\omega''\,\rbrack_{\hat{A}},\qquad \omega,\omega',\omega''\in \hat{X}.\] This will be an algebra: because of the right linearity of $\lbrack\,\cdot\,,\,\cdot\,\rbrack_{\hat{A}}$ we have $b\cdot\lbrack\, \omega,\omega'\,\rbrack_B=\lbrack\, b\cdot \omega,\omega'\,\rbrack_B$. This also makes
the maps of $B$ on $\hat{X}$ commute with the action of
$\hat{A}$. Moreover, if $\omega_1\in \hat{A}$ then \[\lbrack\,  \omega\cdot \omega_1,\omega'\,\rbrack_B   =  \lbrack \,\omega,\omega_1\cdot \omega'\,\rbrack_B,\] which follows from the left $\hat{A}$-linearity of $\lbrack\, \cdot\,,\,\cdot\,\rbrack_{\hat{A}}$. This provides a canonical map $\pi$ from $\hat{X}\underset{\hat{A}}{\otimes} \hat{X}$ to $B$.

\begin{Lem}\label{lem4} The map $S_B:B\rightarrow B: \lbrack\, \omega,\omega'\,\rbrack_B \rightarrow \lbrack\, \theta_X(\omega'),\omega\,\rbrack_B$ is a well-defined bijection.\end{Lem}

\begin{proof} Choose arbitrary $x,y\in X$, and suppose  $\lbrack\, \omega,\omega'\,\rbrack_B=0$. Then \begin{eqnarray*} (\lbrack\,\theta_X(\omega'),\omega\,\rbrack_B\cdot \varphi_X(\,\cdot\, x))(y) &=& \omega'(\theta_X(y_{(0)}))\omega(x_{(0)})\varphi(y_{(1)}x_{(1)})\\ &=& \omega(x_{(0)})\omega'(\sigma_X(y)_{(0)})\varphi (x_{(1)}(\sigma_X(y))_{(1)})\\ &=& (\lbrack\, \omega,\omega'\,\rbrack_B \cdot\varphi_X(y\,\cdot\,))(x)\\ &=& 0,\end{eqnarray*} hence $\lbrack \, \theta_X(\omega'),\omega\,\rbrack_B=0$.
\end{proof}

\noindent This allows us to define a right action of $B$ on $\hat{X}$ by setting $\omega\cdot b=S_B^{-1}(b)\cdot \omega$, since by Lemma \ref{lem5} we have $\omega''\cdot \lbrack\, \omega,\omega'\,\rbrack_B = \lbrack\, \omega'',\omega\,\rbrack_{\hat{A}} \cdot \omega'$.
\begin{Cor} The maps $\lbrack \,\cdot\,,\,\cdot\,\rbrack_B$ and $\lbrack\,\cdot \,,\,\cdot\,\rbrack_{\hat{A}}$ make $\hat{X}$ into a Morita context between $B$ and $\hat{A}$. \end{Cor}

\noindent Remark that this can be used to show that $B$ is a non-degenerate algebra.

\begin{Lem}  For any finite collection $\omega^i\in \hat{X}$ there exists $b\in B$ with $b\cdot\omega^i= \omega^i$.\end{Lem}

\begin{proof}  Put $\omega^i=\varphi_X(\,\cdot \,y_i)$. Then we have to prove that there exist $\omega'^{\,j}$ and $\omega''^{\,j}\in \hat{X}$ such that for any $x\in X$ \[\sum_j\omega'^{\,j}(x_{(0)})\omega''^{\,j}(y_{i(0)})\varphi(x_{(1)}y_{i(1)}) = \varphi_X(xy_i).\] Choose $z\in X$ with $\varphi_X(z)=1$. Put $\omega=\varphi_X(\,\cdot\, z)$ and choose $\omega_1\in \hat{A}$ such that $\omega_1(y_{i(1)})y_{i(0)}z\otimes y_{i(2)} = y_{i(0)}z\otimes y_{i(1)}$ for all $i$. Put $\sum_j \omega'^{\,j}\otimes \omega''^{\,j} = \omega^{(1)}\otimes (\omega^{(2)}\cdot \omega_1)$. Then we have \begin{eqnarray*} \sum_j\omega'^{\,j}(x_{(0)})\omega''^{\,j}(y_{i(0)})\varphi(x_{(1)}y_{i(1)}) &=& \omega(x_{(0)}y_{i(0)}) \varphi(x_{(1)}y_{i(2)})\omega_1(y_{i(1)})\\ &=& \omega(x_{(0)}y_{i(0)}) \varphi(x_{(1)}y_{i(1)})\\ &=& \varphi_X(xy_i).\end{eqnarray*}\end{proof}

\begin{Prop} The projection $\pi:\hat{X}\underset{\hat{A}}{\otimes} \hat{X} \rightarrow B$ is bijective.\end{Prop}

\begin{proof}

\noindent Suppose $\pi(\sum_i \omega^i\underset{\hat{A}}{\otimes} \omega'^i)=0$. Let $b=\sum_j \lbrack\, \omega''^j,\omega'''^{\,j}\,\rbrack$ be a local unit for the $\omega^i$. Then \begin{eqnarray*} \sum_i \omega^i\underset{\hat{A}}{\otimes} \omega'^i &=& \sum_i b\cdot\omega^i\underset{\hat{A}}{\otimes} \omega'^i\\ &=& \sum_{i,j} (\omega''^j\cdot \lbrack\,\omega'''^{\,j},\omega^i\,\rbrack_{\hat{A}})\underset{\hat{A}}{\otimes} \omega'^{\,i}\\ &=& \sum_{i,j} \omega''^j \underset{\hat{A}}{\otimes} (\lbrack\,\omega'''^{\,j},\omega^i\,\rbrack_{\hat{A}} \cdot \omega'^{\,i})\\ &=& \sum_{i,j} \omega''^j \underset{\hat{A}}{\otimes} (\omega'''^{\,j}\cdot \lbrack\,\omega^i, \omega'^{\,i}\,\rbrack_{\hat{B}})\\ &=& 0.\end{eqnarray*}

\end{proof}

\noindent In the following, we will identify $\hat{X}\underset{\hat{A}}{\otimes} \hat{X}$ with $B$.\\

\noindent We can define a right action of $B$ on $X$ by \[x\cdot \lbrack\, \omega,\omega'\,\rbrack_B= \omega(x_{(0)})\omega'(x_{(1)}^{\;\;\;[1]})x_{(1)}^{\;\;\;[2]},\] and then $(b\cdot\omega)(x)=\omega(x\cdot b)$ for all $b\in B$. We can then construct a map \[V_B^t: B\otimes \hat{X}\rightarrow \hat{X}\otimes \hat{X}: b\otimes \omega\rightarrow b\cdot\omega^{(1)}\otimes \omega^{(2)},\] and it will be a bijection. Its inverse is given by \[(V_B^t)^{-1}(\omega\otimes \omega')=\lbrack \, \omega,\omega'^{(1)}\,\rbrack_B \otimes \omega'^{(2)},\] where the first factor is covered by $\hat{A}$-balancedness and a local unit argument.
Note that a same kind of expression can be used for the inverse of $V^t$: we have \[(V^t)^{-1}(\omega\otimes \omega')= \omega^{(1)}\otimes \lbrack\, \omega^{(2)},\omega'\,\rbrack_{\hat{A}}.\]\\

\noindent Now we define a comultiplication on $B$. We first provide some formal intuition. Note that on $\hat{X}$, we have \[\Delta_X(\omega\cdot \omega_1)=\Delta_X(\omega)\cdot \Delta(\omega_1).\] Also, using $(ab)^{[1]}\otimes (ab)^{[2]}=b^{[1]}a^{[1]}\otimes a^{[2]}b^{[2]}$ for $a,b\in A$, we have \[\Delta(\lbrack \,\omega,\omega'\,\rbrack_{\hat{A}}) = \lbrack\, \omega^{(2)},\omega'^{(1)}\,\rbrack_{\hat{A}}\otimes \lbrack\, \omega^{(1)},\omega'^{(2)}\,\rbrack_{\hat{A}}.\] If we then want $\Delta_B(b)\cdot \Delta_X(\omega) = \Delta_X(b\cdot\omega)$, we are led to \[\Delta_B(\lbrack\, \omega,\omega'\,\rbrack_B)= \lbrack \,\omega^{(1)},\omega'^{(2)}\,\rbrack_B\otimes \lbrack\, \omega^{(2)},\omega'^{(1)}\,\rbrack_B.\] We now show that this is a well-defined comultiplication.\\

\noindent Remark that $m=\lbrack\, \omega^{(1)},\omega'^{(2)}\,\rbrack_B\otimes \lbrack\, \omega^{(2)},\omega'^{(1)}\,\rbrack_B$ automatically has a well-defined meaning as a multiplier, using $B$-linearity of $\lbrack\,\cdot\,,\,\cdot\,\rbrack_B$. In fact, for $c\in B$ any of the expressions $(1\otimes c)m, (c\otimes 1)m, m(1\otimes c)$ and $m(c\otimes 1)$ are elements of $B\otimes B$. Since the action of $B$ is unital, we can let $m$ act on the left of $\hat{X}\otimes \hat{X}$. \begin{Lem} For all $\omega'',\omega'''\in \hat{X}$, we have \[m\cdot (\omega''\otimes \omega''')=V^t(b\otimes 1)(V^t)^{-1}(\omega''\otimes \omega'''),\] where $b=\lbrack\, \omega,\omega'\,\rbrack_B$.\end{Lem}  \begin{proof} We have to show that for all $\omega''\in \hat{X}$ and $\omega_1\in \hat{A}$, we have \[V^t(b\cdot \omega''\otimes \omega_1)= m\cdot (\omega''^{(1)}\otimes \omega''^{(2)}\cdot \omega_1).\] Choose $c\in B$. Then \begin{eqnarray*}  ((c\otimes 1)m)\cdot (\omega''^{(1)}\otimes \omega''^{(2)}\cdot \omega_1) &=& \lbrack\, c\cdot \omega^{(1)},\omega'^{(2)}\,\rbrack_B \cdot \omega''^{(1)}\otimes \lbrack \, \omega^{(2)},\omega'^{(1)}\,\rbrack_B\cdot (\omega''^{(2)}\cdot \omega_1)\\ &=& (c\cdot \omega^{(1)})\cdot \lbrack \, \omega'^{(2)},\omega''^{(1)}\,\rbrack_{\hat{A}}\otimes \omega^{(2)}\cdot \lbrack\,\omega'^{(1)},\omega''^{(2)}\cdot \omega_1\,\rbrack_{\hat{A}}\\ &=& (c\cdot \omega^{(1)}\otimes \omega^{(2)})\cdot (\hat{\Delta}_{\hat{A}}(\lbrack\, \omega',\omega''\,\rbrack_{\hat{A}})(1\otimes \omega_1))\\ &=& (c\otimes 1)\cdot (\Delta_X(\omega\cdot \lbrack\, \omega',\omega''\,\rbrack_{\hat{A}})\cdot (1\otimes \omega_1))\\ &=& (c\otimes 1)\cdot (V^t(b\cdot \omega''\otimes \omega_1)).\end{eqnarray*} As $c$ was arbitrary, the lemma is proved.
\end{proof}
\noindent It is then immediate that if we define \[\Delta_B(\lbrack\, \omega,\omega'\,\rbrack_B)=\lbrack\, \omega^{(1)},\omega'^{(2)}\,\rbrack_B\otimes \lbrack\, \omega^{(2)},\omega'^{(1)}\,\rbrack_B,\] then $\Delta_B$ is a well-defined multiplicative comultiplication $B\rightarrow M(B\otimes B)$.\\

\noindent Now we show that  $(B,\Delta_B)$ is an algebraic quantum group. We will do this by explicitly constructing its counit, its antipode and its left invariant functional. This is sufficient by Proposition 2.9. in \cite{VD1}. Define by $\varepsilon_B$ the map \[\varepsilon_B:B\rightarrow k: \lbrack\, \omega,\omega'\,\rbrack_B\rightarrow \omega(1)\omega'(1).\]  \begin{Lem} The functional $\varepsilon_B$ is well-defined, and satisfies the counit property with respect to $\Delta_B$.\end{Lem} \begin{proof} The well-definedness is immediate. Choose $c\in B$, $\omega,\omega'\in \hat{X}$. Then \begin{eqnarray*} (\iota\otimes \varepsilon_B)((c\otimes 1)\Delta_B(\lbrack\, \omega,\omega'\,\rbrack_B))&=& \omega^{(2)}(1)\,\omega'^{(1)}(1)\,\lbrack\, c\cdot \omega^{(1)},\omega'^{(2)}\,\rbrack_B\\ &=& c\cdot \lbrack \,\omega,\omega'\,\rbrack_B.\end{eqnarray*} The other half of the counit property is proven similarly.
\end{proof}

 \begin{Lem} The map $S_B$ defined in Lemma \ref{lem4} satisfies the antipode property.\end{Lem} \begin{proof} For one half of it, we have to prove the following identity: for any $\omega,\omega'\in \hat{X}$ and $c\in B$, \[\lbrack\, \theta_X(\omega'^{(2)}),\omega^{(1)}\,\rbrack_B \cdot \lbrack\, \omega^{(2)},\omega'^{(1)}\cdot c\,\rbrack_B =\omega(1)\omega'(1)c.\] But the left hand expression equals $\lbrack \,\lbrack\,\theta_X(\omega'^{(2)}), \omega^{(1)}\,\rbrack_B\cdot \omega^{(2)},\omega'^{(1)}\cdot c\,\rbrack_B$. Now for $\omega''\in \hat{X}$ we have \[\lbrack\, \omega'', \omega^{(1)}\,\rbrack_B\cdot \omega^{(2)} = \omega(1) \omega'',\] by using Lemma \ref{lem5} and formula $xiv)$ of Proposition \ref{prop2}, so the first expression reduces to $\omega(1)\lbrack\, \theta_X(\omega'^{(2)}),\omega'^{(1)}\cdot c\,\rbrack_B.$ Using formula $xv)$ of proposition \ref{prop2}, we find that for $\omega''\in \hat{X}$ \[\lbrack\,\theta_X(\omega'^{(2)}),\omega'^{(1)}\,\rbrack_B\cdot\omega''= \omega'(1)\omega'',\] proving the identity. The other half of the antipode property follows similarly. \end{proof}

\noindent Now we show that $B$ has a non-zero left invariant functional. Denote by $\;\widehat{}\;$ the map $\hat{X}\rightarrow X: \psi_X(\,\cdot\, x)\rightarrow x$.
 \begin{Lem} The functional \[\varphi_B: B\rightarrow k: \lbrack\, \omega',\omega\,\rbrack_B\rightarrow \omega'(\widehat{\omega})\] is well-defined and left invariant.\end{Lem} \begin{proof}

 If $x\in X$ and $\omega_1\in \hat{A}$, then $\psi_X(\,\cdot \, x)\cdot \omega_1 = \omega_1(S(x_{(1)}))\psi_X(\,\cdot\,x_{(0)})$. Hence if $\omega=\psi_X(\,\cdot\,x)$, then $\widehat{\omega_1\cdot \omega}= \omega_1\cdot\widehat{\omega}$. If also $\omega'\in \hat{X}$, then $(\omega'\cdot \omega_1)(\widehat{\omega}) = \omega'(\widehat{\omega_1\cdot \omega})$. This proves that $\varphi_B$ is well-defined on $B= \hat{X}\underset{\hat{A}}{\otimes} \hat{X}$. \\

 \noindent We prove that $\varphi_B$ is left invariant. First note that the expression \[(\iota\otimes \varphi_B)\Delta_B(\lbrack \,\omega',\omega\,\rbrack_B)\] makes sense as a multiplier of $B$. As the action of $B$ on $X$ is unital, also the expression \[(((\iota\otimes \varphi_B)\Delta_B(\lbrack \,\omega',\omega\,\rbrack_B))\cdot\omega'')(x)\] is meaningful for $\omega''\in \hat{X}$ and $x\in X$. Since the action of $M(B)$ is faithful, it is enough to prove that this expression equals $\omega'(\widehat{\omega})\omega''(x)$. Define \[\xi_{\omega'',x}:B\rightarrow k: b\rightarrow (b\cdot \omega'')(x).\] By unitality we have that $(\xi_{\omega'',x}\otimes \iota)\Delta_B(b)\in B$ for any $b\in B$, and a small calculation yields that for $b=\lbrack \,\omega',\omega\,\rbrack_B$ we have \[(\xi_{\omega'',x}\otimes \iota)\Delta_B(b)= \omega''(x_{(1)}^{\;\;\;\;[2]})\lbrack\, \omega'(x_{(0)}\,\cdot\,),\omega(\,\cdot\,x_{(1)}^{\;\;\;\;[1]})\,\rbrack_B.\] If we then apply $\varphi_B$ and write $\omega=\psi_X(\,\cdot\,y)$, we obtain $\omega''(x_{(1)}^{\;\;\;\;[2]})\omega'(x_{(0)}x_{(1)}^{\;\;\;\;[1]}y)$, which reduces to $\omega'(\widehat{\omega})\omega''(x)$. \end{proof}

\noindent We have proven \begin{Theorem} Together with the map $\Delta_B$ the algebra $B$ is an algebraic quantum group.\end{Theorem}

\subsection*{\qquad \textit{$X$ as a bi-Galois object}}

\noindent Denote by $C$ the vector space $\hat{X}\underset{\hat{A}}{\otimes} X$. Then \[B\times C\rightarrow k : (b ,\omega\otimes x)\rightarrow (b\cdot\omega)(x)\] is a well-defined pairing. Now since the map $\hat{X}\rightarrow X: \psi_X(\,\cdot\, x)\rightarrow x$ is left $\hat{A}$-linear, we can define a natural bijection $B\rightarrow C$ by sending $\lbrack\, \omega,\psi_X(\,\cdot\, x)\,\rbrack_B$ to $\omega\otimes x$. We also have a canonical bijection $\hat{B}\rightarrow B: \varphi_B(\,\cdot \,b)\rightarrow b$. We can use this to identify $C$ with $\hat{B}$. In the following we will denote $\omega\underset{\hat{A}}{\otimes} x$ by $\lbrack\, \omega,x\,\rbrack_C$. Then multiplication in $C$ is essentially defined by \[\lbrack \,\omega^{(1)},x\,\rbrack_C \cdot \lbrack\, \omega^{(2)},y\,\rbrack_C = \lbrack\, \omega,xy\,\rbrack_C.\] Also its counit, antipode and right invariant functional are easily seen to be defined, using the results on the structure of the dual of an algebraic quantum group, by respectively \[\left.\begin{array}{lcl} \varepsilon_C(\lbrack\, \omega,x\,\rbrack_C)&=& \omega(x)\\ S_C(\lbrack\, \varphi_X(\,\cdot\, w),x\,\rbrack_C) &=& \lbrack\, \varphi_X(x\,\cdot\,),w\,\rbrack_C, \\ \psi_C(\lbrack \,\omega,x\,\rbrack_C)&=&\omega(1)\psi_X(x).\end{array}\right.\] There does not seem to be any nice form for the comultiplication on $C$.\\

\noindent It is immediate that the right action of $B$ on $X$ makes $X$ a unital right $B$-module algebra. Then we know from \cite{VDZ1} that there is a left coaction of $C$ on $X$. If we denote it by $\gamma$, then we have the formula \[(\iota\otimes \omega)\gamma(x)=\lbrack \omega,x\rbrack_C.\] It is also clear that this coaction makes $X$ a left $C$-Galois object, since the adjoint of the Galois map is exactly $V_B^t$. As $X$ is a $\hat{A}$-$B$-bi-module, it is clear that the
coactions of $C$ and $A$ commute. Hence

\begin{Theorem} The algebra $X$ is a
$C$-$A$-bi-Galois object, i.e.\! simultaneously a left $C$-Galois object and right $A$-Galois object such that the coactions of $C$ and $A$ commute. The map $\varphi_X$ will be invariant for the coaction of $C$, while $\psi_X$ will be $\delta_C^{-1}$-invariant.\end{Theorem}
\begin{proof} We have already shown the validity of the first statement. As for the invariance of $\varphi_X$, we have to show that $\varphi_X(x\,\cdot\, b)=\varphi_X(x)\varepsilon_B(b)$ for $x\in X$ and $b\in B$. Choosing $\omega\in \hat{X},y\in X$, we have for all $x\in X$ that \begin{eqnarray*} \varphi_X(x\,\cdot\, \lbrack \, \omega,\varphi_X(y\,\cdot\,)\,\rbrack_B) &=& \omega( x_{(0)})\varphi(y_{(1)}S(x_{(1)}))\varphi_X(y_{(0)})\\  &=& \varphi_X(y)\varphi_X(x)\omega(1)\\ &=& \varphi_X(x)\varepsilon(\lbrack\, \omega,\varphi_X( y\cdot)\,\rbrack_B).\end{eqnarray*}

\noindent As for the $\delta_C^{-1}$-invariance of $\psi_X$, this follows from the fact that \[\psi_X(x)1=(\psi_C\otimes \iota)(\gamma(x)).\] This shows that $\psi_X$ bears the same relation to $\psi_C$ as $\varphi_X$ did to $\varphi$, and reasoning by duality the claim follows.\end{proof}

\noindent We give another characterization of the algebra $B$. When $x\in X$, we will denote by $R_x$ the map `right multiplication with $x$'. By $B_0(X)$ we denote $X\otimes \hat{X}$ seen as finite rank operators on $X$ in the natural way.

 \begin{Prop} The algebra $B$ consists exactly of those maps $F:X\rightarrow X$ which commute with the left action of $\hat{A}$ and such that $\{R_xF\mid x\in X\}\subseteq B_0(X)$\end{Prop} \begin{proof} It is not difficult to show that every element of $B$ satisfies this condition. For the other way, we first show that $F$ can be seen as an element of $M(B)$. The commutation of $F$ with the left action of $\hat{A}$ let's us identify $F$ with a functional $\omega_F$ on $C= \hat{X}\underset{\hat{A}}{\otimes} X$ by sending $\omega\otimes x$ to $\omega(F(x))$. As such, for any $c\in C$, we have $(\omega_F\otimes \iota)(\Delta_C(c))\in C$ and $(\iota\otimes \omega_F)(\Delta_C(c))\in C$, for if $c= \langle \varphi_X(\,\cdot\, y),x\rangle_C$, then \begin{eqnarray*} (\omega_F\otimes \iota)(\Delta_C(c))&=& \omega_F(x_{(-2)})\varphi_X(x_{(0)}y) x_{(-1)} \\ &=& \omega_F(x_{(-1)})\varphi_X(x_{(0)}y_{(0)})S_C^{-1}(y_{(-1)})\\ &=& \varphi_X(F(x)y_{(0)})S_C^{-1}(y_{(-1)})\in C,\end{eqnarray*} while \begin{eqnarray*} (\iota\otimes \omega_F)(\Delta_C(c))&=& \varphi_X(x_{(0)}y)\omega_F(x_{(-1)}) x_{(-2)}\\ &=& \varphi_X((R_yF)(x_{(0)})) x_{(-1)}\in C.\end{eqnarray*} The remark before proposition 4.3. \!\!of \cite{VD1} let's us conclude that $\omega_F\in M(B)$, and hence $F$ is the right action by $\omega_F$. As the map sending $m\otimes x \in M(B)\otimes X$ to $R_xm$ is seen to be injective, and as $B\otimes X\rightarrow B_0(X): b\otimes x\rightarrow R_xb$ is seen to be a bijection, we conclude that $\omega_F\in B$.\end{proof}

 \begin{Cor} If $C'$ is another algebraic quantum group making $X$ a $C'-A$-bi-Galois object, then $C$ and $C'$ are isomorphic as algebraic quantum groups. \end{Cor}
 \begin{proof} It is enough to proof that the dual $B'$ of $C'$ is isomorphic with $B$. But the previous proposition implies that $\pi(B')\subseteq B$, with $\pi$ the associated faithful representation of $B'$ as operators on $\hat{X}$. Since the Galois property forces the natural map $B'\otimes \hat{X}\rightarrow \hat{X}\otimes \hat{X}$ to be an isomorphism, we have $\pi(B')=B$.
 \end{proof}


\subsection*{\qquad \textit{Concerning $^*$-structures}}

\noindent Suppose now again that $A$ is $^*$-algebraic (over
$\mathbb{C}$), $X$ a $^*$-algebra and $\alpha$ a $^*$-homomorphism. Then we know that $\varphi_X$ is a positive
functional on $X$. We can introduce on $\hat{X}$ the
$^*$-operation \[\omega^*(x)=\overline{\omega(x^*)}.\]  \begin{Lem}\label{lem6} For all $\omega\in \hat{X}$ and $\omega_1\in \hat{A}$, we have \[(\omega\cdot \omega_1)^* = \omega_1^*\cdot \omega^*.\]\end{Lem} \begin{proof} Choose $x\in X$, then \begin{eqnarray*} (\omega\cdot \omega_1)^*(x) &=& \overline{(\omega\cdot \omega_1)(x^*)}\\ &=& \overline{\omega(x_{(0)}^*)\omega_1(x_{(1)}^*)}\\ &=& \omega^*(x_{(0)})(\hat{S}^{-1}(\omega_1^*))(x_{(1)})\\ &=& (\omega_1^*\cdot \omega^*)(x).\end{eqnarray*}\end{proof}

\noindent Now consider
the following $^*$-operation on the dual $C$ of $B$: \[\lbrack\, \omega,x\,\rbrack_C^*:=
\lbrack\, \omega^*,x^*\,\rbrack_C.\] By the previous lemma, it is seen to be a
well-defined involution. Writing elements of $C$ in the form
$(\iota\otimes \omega)(\gamma(x))$, and using the formula
\[(\omega^*)^{(1)}\otimes (\omega^*)^{(2)}= (\omega^{(2)})^*\otimes
(\omega^{(1)})^*,\] we find that $^*$ is
anti-multiplicative on $C$. Finally, $\Delta_C$ is $^*$-preserving,
which is again easily seen by writing an element of $C$ in the
aforementioned form and using that $\gamma$ is $^*$-preserving.\\

\noindent Now we show that $\psi_C$ is positive. Take $c\in C$ and
$z\in X$ with $ \varphi_X(z^*z)=1$. Write $c\otimes z$ as $\sum_i
\gamma(x_i)(1\otimes y_i)$, then $c^*c= \sum_{i,j} (\iota\otimes
\varphi_X((y_i)^*\cdot y_j))(\gamma((x_i)^*x_j))$, so \[c^*c =
\sum_{i,j}\langle \varphi_X((y_i)^*\cdot y_j),(x_i)^*x_j\rangle_C.\]
Applying $\psi_C$, we get \begin{eqnarray*} \psi_C(c^*c) &=&
\sum_{i,j} \varphi_X((y_i)^*y_j)\psi_X((x_i)^*x_j)\\ &\geq&
0,\end{eqnarray*} since the matrices
$(a_{i,j})=(\varphi_X((y_i)^*y_j))$ and
$(b_{i,j})=(\psi_X((x_i)^*x_j))$ are both positive definite. Hence

\begin{Theorem} If $A$ is a $^*$-algebraic quantum group and $(X,\alpha)$ is a Galois object such that $X$ is a $^*$-algebra and $\alpha$ is $^*$-preserving, then the algebraic quantum group $(C,\Delta_C)$ is a $^*$-algebraic quantum group.\end{Theorem}

\noindent We can also consider the corresponding $^*$-operation on $B$. Then we have a formula corresponding to the one of the above lemma, namely for $\omega\in \hat{X}$ and $b\in B$, we have \[ (b\cdot \omega)^*= \omega^*\cdot b^*.\]

\subsection*{\qquad \textit{Monoidal equivalence}}

\noindent \begin{Def} If $D$ is an algebraic quantum group, we denote by $D-Rep$ the monoidal category of unital left $D$-modules, where a morphism between two objects $V$ and $W$ is a linear map $V\rightarrow W$ which intertwines the module structure. \end{Def} \noindent Note that this category is indeed monoidal by the usual tensor product of two modules, because of the unitality assumption.\\

\noindent We construct a monoidal equivalence between $\hat{A}-Rep$ and $B-Rep$. This equivalence is given by the
functor $F=\hat{X}\underset{\hat{A}}{\otimes}-$.
 This functor is
monoidal with respect to the natural isomorphism
\[n_{\otimes}:\omega\underset{\hat{A}}{\otimes}(v\otimes w)\rightarrow
(\omega^{(1)}\underset{\hat{A}}{\otimes} v)\otimes
(\omega^{(2)}\underset{\hat{A}}{\otimes} w),\] where it is clear how to
interpret this and how to show bijectivity. \\

\noindent We briefly argue how to construct the inverse, without entering into details. Remark that $F$ sends the left $\hat{A}$-module algebra $A$  to $\hat{X}\underset{\hat{A}}{\otimes} A$, which is naturally identified with the algebra $X^{\textrm{op}}$ via $\omega\otimes a\rightarrow \omega(a^{[2]})a^{[1]}$. The induced left $B$-module structure on $X^{\textrm{op}}$ is dual to the right one on $\widehat{X^{\textrm{op}}}=\hat{X}$. This then makes $X^{\textrm{op}}$ a right $C$-Galois object, and we can apply on $B-Rep$ the functor $ \widehat{X^{\textrm{op}}}\underset{B}{\otimes}-$. Then $\widehat{X^{\textrm{op}}}\underset{B}{\otimes}( \hat{X}\underset{\hat{A}}{\otimes}-)$ is seen to be naturally equivalent with the identity, using first the isomorphism $\widehat{X^{\textrm{op}}}\underset{B}{\otimes} \hat{X}\rightarrow \hat{A}: \omega\otimes \omega'\rightarrow \lbrack\, \omega,\omega'\,\rbrack_{\hat{A}}$, then by using that $\hat{A}\underset{\hat{A}}{\otimes}-$ is simply the identity, using a local unit argument.\\

\noindent Suppose now again that $A$ is $^*$-algebraic (over
$\mathbb{C}$), $X$ a $^*$-algebra and the coaction a $^*$-homomorphism. We
can equip the category of unital left $\hat{A}$-modules with an anti-linear
conjugation functor $Conj_{\hat{A}}$ by sending
$(V,\pi)$ to $(\overline{V},\overline{\pi})$, where $\overline{V}$ is the conjugate vector space of $V$ and with
\[\overline{\pi}(\omega_1) \cdot \overline{v} =
\overline{\hat{S}(\omega_1^*)\cdot v}.\] There is a natural
isomorphism \[n_{Conj}:(Conj_B \circ
(\hat{X}\underset{\hat{A}}{\otimes} -) \circ Conj_{\hat{A}})\rightarrow
\hat{X}\underset{\hat{A}}{\otimes} -,\] given by
\[\overline{\omega\underset{\hat{A}}{\otimes} \overline{v}}\quad \rightarrow \quad \theta_X(\omega^*) \underset{\hat{A}}{\otimes} v.\] To see that this is well-defined, we have to prove that \[\theta_X((\omega\cdot \omega_1)^*)\underset{\hat{A}}{\otimes} v = \theta_X(\omega^*)\underset{\hat{A}}\otimes \hat{S}(\omega_1^*)\cdot v,\] for all $\omega\in \hat{X},\omega_1\in \hat{A}$ and $v\in V$. Now $(\omega\cdot \omega_1)^*= \omega^*\cdot \hat{S}^{-1}(\omega_1^*)$ by Lemma \ref{lem6}, and $\theta_X(\omega\cdot \omega_1)=\theta_X(\omega)\cdot \hat{S}^2(\omega_1)$ by $vii)$ of Proposition \ref{prop2}, which proves the identity. To prove that it is a natural transformation, we have to show that \[\theta_X((S_B(b^*)\cdot \omega)^*)= b\cdot \theta_X(\omega^*)\] for $b\in B$ and $\omega\in \hat{X}$. But $\theta_X((S_B(b^*)\cdot \omega)^*)=\theta_X(S_B^{-2}(b)\cdot \omega^*)$, so the identity follows by the `mirror version' of the identity $vii)$ in Proposition \ref{prop2}.\\

\noindent Now look at the category $\hat{A}-Rep^*$ of unital left $\hat{A}$-modules which
have a pre-Hilbertspace structure such that the resulting representation of $\hat{A}$ is $^*$-preserving. The morphisms in the category are now required to have an adjoint. If $V$ is an object in this category, we have a canonical morphism
$\textrm{Sc}:V\otimes \overline{V}\rightarrow \varepsilon$, sending $v\otimes \overline{w}$ to $\langle v,w\rangle$. Note that we take the convention where the scalar product is antilinear on its second argument. Using the natural
isomorphisms for tensoring and conjugating, the map Sc is sent to a map $(\hat{X}\underset{\hat{A}}{\otimes} V)\otimes
\overline{(\hat{X}\underset{\hat{A}}{\otimes} V)}$. This provides an sesquilinear form on $\hat{X}\underset{\hat{A}}{\otimes} V$. We will show in Proposition \ref{prop3} that it equips $\hat{X}\underset{\hat{A}}{\otimes} V$ with a pre-Hilbertspace structure. This is the \textit{categorical approach} to arrive at the induced Hilbert space structure. \\

\noindent We can also use a specific $\hat{A}$-valued inner product on $\hat{X}$ to perform the induction, given by \[\langle \omega,\omega'\rangle_{\hat{A}} = \lbrack\, \omega^*,\omega'\,\rbrack_{\hat{A}}.\] The fact that $\langle\, \cdot\,,\,\cdot\,\rangle_{\hat{A}}$ is a $\hat{A}$-valued inner product means that it is a sesquilinear map with values in $\hat{A}$, antilinear in the \textit{first} argument, such that \begin{enumerate} \item $\langle \omega,\omega'\rangle_{\hat{A}}^*=\langle \omega',\omega\rangle_{\hat{A}}$, \item $\langle \omega,\omega\rangle_{\hat{A}}\geq 0$ with equality iff $\omega=0$, and \item $\langle \omega,\omega'\cdot \omega_1\rangle_{\hat{A}}=\langle \omega,\omega'\rangle_{\hat{A}}\cdot \omega_1$ for $\omega_1\in \hat{A}$.\end{enumerate} Moreover, we have \begin{enumerate}\item[4.] $\langle b\cdot \omega,\omega'\rangle_{\hat{A}} = \langle \omega,b^*\cdot \omega'\rangle_{\hat{A}}$, $b\in B$.\end{enumerate} We will prove this below. We can then make an inner product on $\hat{X}\underset{\hat{A}}{\otimes} V$  by setting \[\langle \omega\otimes v,\omega'\otimes w\rangle = \langle v,\langle \omega,\omega'\rangle_{\hat{A}}\cdot w\rangle.\] This is the \textit{C$^*$-algebraic} approach (cf. \cite{Rie1}).\\

\noindent Finally, there is a \textit{von Neumann algebraic} approach. Namely, introduce in $\hat{A}$ the inner product determined by $\langle \omega_1,\omega_2\rangle = \hat{\psi}(\omega_2^*\omega_1)$ (i.e.\! $\langle \varphi(\,\cdot\, a),\varphi(\,\cdot\, b)\rangle = \varphi(b^*a)$), and introduce in $\hat{X}$ the inner product determined by $\langle \varphi_X(\,\cdot\, y),\varphi_X(\,\cdot\, x)\rangle = \varphi_X(x^*y)$. For $\omega\in\hat{X}$, denote by $L_\omega$ the map $\hat{A}\rightarrow \hat{X}$ which sends $\omega_1$ to $\omega\cdot \omega_1$. Then $L_\omega^*L_{\omega'}$ will be left multiplication with some element of $\hat{A}$, and identifying the operator with this element, we can define \[\langle \omega\otimes v,\omega'\otimes w\rangle = \langle v,L_\omega^*L_{\omega'}\cdot w\rangle.\] Remark however that with this scalar product on $\hat{X}$, the right or left representation of $\hat{A}$ on $\hat{X}$ is in general \textit{not} a $^*$-representation. This is because the left action of $\hat{A}$ on $\hat{X}$ is not an analogue of the left multiplication of $\hat{A}$ on $\hat{A}$.

\begin{Prop}\label{prop3} All three sesquilinear forms on $\hat{X}\underset{\hat{A}}{\otimes} V$ coincide, providing this space with a pre-Hilbertspace structure such that the induced representation of $B$ is $^*$-preserving.\end{Prop}

\begin{proof}
\noindent We first show that $L_\omega^*L_{\omega'} = \langle \omega,\omega'\rangle_{\hat{A}}$. This means that for all $\omega,\omega'\in \hat{X}$ and $\omega_1,\omega_2\in \hat{A}$ we have \[\langle \omega\cdot \omega_1,\omega'\cdot \omega_2\rangle = \langle \omega_1,\langle\omega,\omega'\rangle_{\hat{A}}\omega_2\rangle.\] Writing $\omega'=\varphi_X(\,\cdot\, x)$ and $\omega_2=\varphi(\,\cdot\, a)$, we have \begin{eqnarray*} \langle \omega\cdot \omega_1,\omega'\cdot \omega_2\rangle &=& \langle \omega\cdot \omega_1,\varphi_X(\cdot x_{(0)})\varphi(\delta S(x_{(1)})a)\rangle
 \\ &=& (\omega\cdot \omega_1)(x_{(0)}^*)\varphi(a^*S(x_{(1)})^*\delta),\end{eqnarray*} while  \begin{eqnarray*}  \langle \omega_1,\langle\omega,\omega'\rangle_{\hat{A}}\omega_2\rangle &=& \langle \omega_1,\omega^*(x_{(0)})\varphi(\cdot x_{(1)})\varphi(\delta S(x_{(2)})a)\rangle
 \\&=& \omega(x_{(0)}^*)\varphi(a^*S(x_{(2)})^*\delta)\omega_1(x_{(1)}^*),\end{eqnarray*} which shows the equality.\\

 \noindent Now we show the equivalence of the categorical and the C$^*$-algebraic approach. In the categorical approach, we have for $\omega\in \hat{X}$ that \begin{eqnarray*} \langle \omega^{(1)}\cdot \omega_1\otimes v,(\omega^{(2)}\circ \theta_X^{-1})^*\otimes w\rangle &=& F\textrm{(Sc)}(\omega\otimes v\otimes  w) \\ &=& \omega(1)\langle v, w\rangle,\end{eqnarray*} where $F$ still denotes the functor $\hat{X}\underset{\hat{A}}{\otimes}-$. Since \[\langle (\theta_X^{-1}(\omega^{(2)}))^* , \omega^{(1)}\rangle_{\hat{A}} = \omega(1)\varepsilon,\] we are done.\\

 \noindent From the von Neumann algebraic picture, it follows that this inner product is indeed positive (using that the $^*$-structure on $\hat{A}$ is exactly the adjoint of left multiplication with respect to the stated scalar product on $\hat{A}$). From the categorical picture, it follows quite immediately that the resulting representation of $B$ is $^*$-preserving. Finally, also the non-degeneracy of the inner product follows from the categorical viewpoint: for the vectors of length zero in $F(V)$ form a subobject of $V$, which is sent to $0$ by $F^{-1}$.\\

\noindent The fact that the $\hat{A}$-valued map is indeed a $\hat{A}$-hermitian product follows then in a straightforward manner from the categorical and von Neumann algebraic viewpoint.

\end{proof}

\noindent So in any case, $\hat{X}\underset{\hat{A}}{\otimes}-$ can be lifted to a functor between the
$^*$-representation categories.

\begin{Prop} The natural transformation $n_\otimes$ is unitary.\end{Prop}
\begin{proof} This follows by using the C$^*$-algebraic picture, and using the identity \[\hat{\Delta}_{\hat{A}}(\langle \omega,\omega'\rangle_{\hat{A}})= \langle \omega^{(1)},\omega'^{(1)}\rangle_{\hat{A}}\otimes \langle \omega^{(2)},\omega'^{(2)}\rangle_{\hat{A}}.\] \end{proof}

\noindent The map $n_{Conj}$ however will only be unitary in case $S^2=\iota$, and moreover, if this is not the case, no unitary intertwiner between $\overline{\hat{X}\otimes \overline{V}}$ and $\hat{X}\otimes V$ can be constructed. We can however repair this situation by changing the definition of the conjugation operator: now we send $(V,\pi)$ to $(\overline{V},\overline{\pi})$ with
\[\overline{\pi}(\omega_1) \cdot \overline{v} =
\overline{\hat{R}(\omega_1)^*\cdot v},\] where $\hat{R}$ is the unitary antipode for $\hat{A}$ (see \cite{PRo1}). Note that we can take the square $\theta_X^{1/2}$ of the positive diagonizable operator $\theta_X$, i.e.\! $\theta_X^{1/2}:X\rightarrow X$ is a diagonizable map $X\rightarrow X$ with positive eigenvalues, and $\theta_X^{1/2}(\theta_X^{1/2}(x))=\theta_X(x)$ for all $x\in X$. This $\theta_X^{1/2}$ will still be a multiplicative automorphism of $X$.
\begin{Prop} The map  \[\overline{\omega\underset{\hat{A}}{\otimes} \overline{v}} \rightarrow \theta_X^{1/2}(\omega^*) \underset{\hat{A}}{\otimes} v\] provides a well-defined unitary intertwiner between $\overline{\hat{X}\otimes \overline{V}}$ and $\hat{X}\otimes V$.\end{Prop}
\begin{proof} \noindent To see if this map is well-defined, we have to check the identity \[\theta_X^{1/2}((\omega\cdot \omega_1)^*)= \theta_X^{1/2}(\omega^*)\cdot \hat{R}(\omega_1^*),\] with $\omega\in \hat{X},\omega_1\in \hat{A}$. Using Lemma \ref{lem6}, this reduces to proving that \[\theta_X^{1/2}(\omega\cdot \omega_1) = \theta_X^{1/2}(\omega)\cdot \theta_A^{1/2}(\omega_1),\] where $\theta_A=\hat{S}^2$ and $\theta_A^{1/2}$ is its positive root (so that $\theta_A^{1/2}\circ \hat{R}=\hat{S}$). This follows from $\theta_X(\omega\cdot \omega_1) = \theta_X(\omega)\cdot \hat{S}^2(\omega_1)$, by taking $\omega$ an eigenvector for $\theta_X$ and $\omega_1$ an eigenvector for $\hat{S}^2$.\\

\noindent The way to prove that this is a natural transformation is very similar, and we omit it.\\

\noindent Finally, to prove unitarity we have to prove the identity \[\hat{R}(\lbrack \, \omega',\omega\,\rbrack_{\hat{A}})=\lbrack \,\theta_X^{-1/2}(\omega),\theta_X^{1/2}(\omega')\,\rbrack_{\hat{A}},\] where $\omega,\omega'\in \hat{X}$. Applying $\theta_A^{1/2}$ and using Lemma \ref{lem6}, this reduces to proving that \[\lbrack\, \theta_X^{1/2}(\omega),\theta_X^{1/2}(\omega')\,\rbrack_{\hat{A}} = \theta_A^{1/2}(\lbrack\, \omega,\omega'\,\rbrack_{\hat{A}}).\] This identity is true when $\theta_X^{1/2}$ is replaced by $\theta_X$ and $\theta_A^{1/2}$ is replaced by $\hat{S}^2$, using formula $xvi)$ of Proposition \ref{prop2}. Again by using an eigenvector argument, it is also true as stated.
\end{proof}

\section{An example}

\noindent The following examples of infinite-dimensional Hopf algebras with a left invariant functional can be found in \cite{VD1} and \cite{VDZ2}. We slightly generalize the construction to fit them both in a family.

 \begin{Def} Let $n>1$, $m\geq 1$ be natural numbers, and $\lambda\in k$ such that $\lambda^m$ is a primitive $n$-th root of unity. Let $A^{n,m}_\lambda$ be the unital algebra over $k$ generated by elements $a$, $a^{-1}$ and $b$ with $a^{-1}$ the inverse of $a$, $ab=\lambda ba$ and $b^n=0$. Then we can define a comultiplication on $A^{n,m}_\lambda$ determined on the generators by \[\Delta(a)=a\otimes a,\]\[\Delta(b)=b\otimes a^m+1\otimes b.\] This makes $A^{n,m}_\lambda$ an algebraic quantum group of compact type.\end{Def} \noindent To proof that this comultiplication is indeed well-defined, we only have to use the well-known fact that $(s+t)^l=s^l+t^l$ when $s,t$ are variables satisfying the commutation $st=qts$ with $q$ a primitive $l$-th root of unity. Now $A^{n,1}_\lambda$ is the example in \cite{VDZ2}, and with the further relation $a^n=1$, this reduces to the two-generator Taft algebras. The Hopf algebra $A^{n,2}_\lambda$ is isomorphic with the example constructed in \cite{VD1}. \\

 \noindent The left invariant functional $\varphi$ of $A=A^{n,m}_\lambda$ is defined by $\varphi(a^pb^q)=\delta_{p,0}\delta_{q,n-1}$, $p\in \mathbb{Z}$, $0\leq q<n$. As $A$ is infinite-dimensional, its dual $\hat{A}$ is necessarily of discrete type and not compact, i.e.\! it is a genuine multiplier Hopf algebra. This is a difference with the Taft algebras, which are self-dual. Remark that there can still be defined a pairing between $A$ and itself, but it will be degenerate.\\


\noindent In \cite{Mas1} the Galois objects for the Taft algebras were classified. It provides the motivation for the following construction. Fix $A=A^{n,m}_\lambda$ as above, and assume moreover that $\lambda$ is a primitive $n$-th root of unity and $m$ and $n$ are coprime. The condition `$\lambda^m$ is a primitive root of unity' follows from this assumption.

 \begin{Def} Take $\mu \in k$. Let $X=X^{n,m}_{\lambda,\mu}$ be the unital algebra generated by $x$, $x^{-1}$ and $y$, with $x^{-1}$ the inverse of $x$, such that $xy=\lambda yx$ and $y^n=\mu x^{mn}$. A right coaction $\alpha$ of $A$ on $X$ is defined on the generators by \[\alpha(x)=x\otimes a,\] \[\alpha(y)= y\otimes a^m+1\otimes b.\]\end{Def} \noindent It is again easy to show that this has a well-defined extension to the whole of $X$.

 \begin{Prop} $(X,\alpha)$ is a right $A$-Galois object.\end{Prop}\begin{proof} First of all, we have to see if $X$ is not trivial. We follow standard procedure. Let $V$ be a vector space over $k$ which has a basis of vectors of the form $e_{p,q}$ with $p\in \mathbb{Z}$ and $0\leq q<n$. Define  operators $x'$ and $y'$ by \[\left.\begin{array}{lcrll} x'\cdot e_{p,q}&=&e_{p+1,q}&&  \textrm{for all } p\in \mathbb{Z},0\leq q<n,\\y'\cdot e_{p,q}&=& \lambda^{-p} \; e_{p,q+1}&& \textrm{if } p\in \mathbb{Z}, 0\leq q<n-1,\\y'\cdot e_{p,n-1}&=&\mu\lambda^{-p} \;e_{p+nm,0}&& \textrm{if } p\in \mathbb{Z}.\end{array}\right.\] Then it is easy to see that $x'$ is invertible and that $x'y'=\lambda y'x'$. Also: \begin{eqnarray*} y'^n \cdot e_{p,q} &=& \lambda^{-p(n-1-q)}y'^{1+q}\cdot e_{p,n-1}\\ &=& \mu \lambda^{-p(n-1-q)}\lambda^{-p}y'^q\cdot e_{p+nm,0}\\ &=& \mu\lambda^{-p(n-1-q)}\lambda^{-p}\lambda^{-pq}e_{p+nm,q}\\ &=& \mu\lambda^{-pn}e_{p+nm,q}\\ &=& \mu x'^{mn}\cdot e_{p,q}.\end{eqnarray*} This gives us a non-trivial representation of $X$. Moreover, it is easy to see that this representation is faithful.\\

 \noindent Define by $\tilde{\beta}:A\rightarrow X^{\textrm{op}}\otimes X$ the homomorphism generated by \[\tilde{\beta}(a)=x^{-1}\otimes x,\]\[\tilde{\beta}(b)=- yx^{-m}\otimes x^m +1\otimes y.\] This is well-defined: for example, we have \begin{eqnarray*} \tilde{\beta}(b)^n &=& (-yx^{-m})^n\otimes x^{mn}+\mu 1\otimes x^{mn}\\ &=& ((-1)^n \lambda^{mn(n-1)/2}+1)\mu(1\otimes x^{mn})\\ &=& 0,\end{eqnarray*} using that $\lambda^m$ is a primitive root of unity.
  Denoting $\beta=(f^{-1}\otimes \iota)\tilde{\beta}$ with $f$ the canonical map $X\rightarrow X^{\textrm{op}}$ and writing $\beta(c)=c^{[1]}\otimes c^{[2]}$ for $c\in A$, it is easy to compute that \[z_{(0)}z_{(1)}^{\;\;\;[1]}\otimes z_{(1)}^{\;\;\;[2]}=1\otimes z,\] \[c^{[1]}c^{[2]}_{\;\;\;(0)}\otimes c^{[2]}_{\;\;\;(1)}=1\otimes c\] for all $z\in \{x,y,x^{-1}\}$ and $c\in \{a,b,a^{-1}\}$, and hence for all $z\in X, c\in A$. This is enough to know that the action is Galois.\end{proof} \noindent The extension $k\subseteq X$ will be cleft (see e.g.\! Definition 2.2.3. in \cite{Sch1}), by the comodule isomorphism $\Psi_X:X\rightarrow A: x^py^q\rightarrow a^pb^q$, $p\in\mathbb{Z}$ and $0\leq q<n$. The associated scalar cocycle $\eta$ is given by $\eta(a^pb^q\otimes a^rb^s)=0$, except for $q=s=0$, where it is 1, and when $q+s=n$, in which case it equals $\mu\lambda^{-rq}$.\\

\noindent We determine the extra structure occurring in this example. First note that we have shown that the elements of the form $x^py^q$ with $p\in \mathbb{Z}$ and $0\leq q< n$ form a basis. Then we have  \[\left.\begin{array}{llllll} \varphi_X(x^py^q)&=&\delta_{q,n-1}\delta_{p,0}&&  \textrm{for } &p\in \mathbb{Z}, 0\leq q< n,\\ \psi_X(x^py^q)&=&\delta_{q,n-1}\delta_{p,m(1-n)} \lambda^{-m} && \textrm{for } & p\in \mathbb{Z}, 0\leq q< n, \\ &&&&\\ \delta_X&=&x^{(n-1)m} &&\\ &\\ \sigma_X(x)&=&\lambda^{-1}x,&\theta_X(x)&=&x \\ \sigma_X(y)&=&y, &  \theta_X(y)&=&\lambda^m y.\end{array}\right.\] It is of course the nature of the example which makes the structure so similar to the one of $A$. However: note that the formula for the antipode in $A$, namely $S(a)=a^{-1}$ and $S(b)=-ba^{-m}$, has no well-defined analogue in $X$: for if we set $S_X(x)=x^{-1}$ and $S_X(y)=-yx^{-m}$, then we would expect $\mu x^{-mn}=S_X(y^n)=S_X(y)^n$, but $S_X(y)^n=(-yx^{-m})^n= (-1)^n\lambda^{mn(n-1)/2} y^n x^{-mn}=-\mu$, which shows that $S_X$ can not be extended to an anti-morphism of $X$ (unless $\mu=0$).\\

\noindent Now we determine the associated algebraic quantum group $C$. Note that we could determine the structure with the help of the cocycle, but we wish to directly use the Galois object itself, since this seems easier. In particular, we exploit the duality between $C$ and $B$. \\

\noindent We first give a heuristic reasoning. We determine the algebra structure of $B$. We need a description of the dual $\hat{A}$ of $A^{n,m}_\lambda$. It has a basis consisting of expressions $e_pd^q$ with $p\in \mathbb{Z}$ and $0\leq q< n$, where $e_p\in \hat{A}$ and $d\in M(\hat{A})$, such that $e_pe_q=\delta_{p,q} e_p$, $de_p=e_{p+m}d$ and $d^n=0$. With $c=\sum_k \lambda^{-k}e_k\in M(\hat{A})$, the comultiplication is determined by \[\Delta(e_p)=\sum_t e_t\otimes e_{p-t},\]\[\Delta(d)=d\otimes c+1\otimes d.\]

\noindent Now the left action of $\hat{A}$ on $X$ is given by \[\left.\begin{array}{lcll} e_s\cdot x^py^q&=&\delta_{p,s-mq} \;x^py^q, \\d\cdot x^py^q&=& C_q \;x^py^{q-1}, &  q>0, \\d\cdot x^p&=&0,\end{array}\right.\] where $C_q$ denotes some binomial coefficient of q-calculus (cf.\! \cite{Kli1}). Consider the operators $g_s$ and $h$ acting on the right of $X$ by \[\left.\begin{array}{rcll} x^py^q\cdot g_s&=&\delta_{p,-s} \; x^py^q, \\x^py^q\cdot h&=&C_q \;x^{p+m}y^{q-1}, & q>0, \\x^p \cdot h &=&0.\end{array}\right.\] Then it is easy to see that $h$ and $g_s$ commute with the left action of $\hat{A}$. We see that $h\cdot g_s= g_{s+m}\cdot h$, that $g_sg_t=\delta_{s,t}g_s$ and that $h^n=0$. The span of $g_sh^q$ will form our algebra $B$. Now denote by $u_{p,q}$ the elements in $C$ such that $\langle u_{p,q},e_rd^s\rangle = \delta_{p,r}\delta_{q,s}$, and denote $u=u_{-1,0}$, $v=u_{0,0}$ and $w=u_{0,1}$. Then we have $\gamma(x)=u\otimes x$ and $\gamma(y)=v\otimes y+w\otimes x^m$ by using the action of $B$. Since this has to commute with $\alpha$, we find that $v=1$. Using that $y^n=\mu x^{mn}$ we find that $\mu+w^n = \mu u^{mn}$, and using $xy=\lambda yx$, we get $uw=\lambda wu$. Furthermore, the fact that $x$ is invertible gives that $u$ is invertible. This then completely determines the structure of $C$. The coalgebra structure is determined by the usual \[\Delta(u)=u\otimes u,\]\[\Delta(w)= w\otimes u^m+1\otimes w.\]

\noindent We can now make things exact. \begin{Prop} Let $C$ be the unital algebra generated by three elements $u,u^{-1}$ and $w$, such that $u^{-1}$ is the inverse of $u$, $uw=\lambda wu$ and $\mu\cdot 1+w^n = \mu u^{mn}$. Then $C$ is not trivial. We can define a unital multiplicative comultiplication $\Delta_C$ on $C$, given on the generators by \[\Delta_C(u)=u\otimes u,\]\[\Delta_C(w)= w\otimes u^m+1\otimes w,\] making it an algebraic quantum group of compact type. It has a left coaction $\gamma$ on $X$ determined by \[\gamma(x)=u\otimes x,\]\[\gamma(y)=1\otimes y+w\otimes x^m,\] making it a $C$-$A$-bi-Galois object.\end{Prop}

\begin{proof} It is easy to see that $\Delta_C$ and $\gamma$ can be extended, that $\Delta_C$ is coassociative and $\gamma$ a coaction, and that $\gamma$ commutes with the right coaction of $A$. Since now $C$ is already a bi-algebra, it follows from the general theory of Hopf-Galois extensions that if $\gamma$ can be shown to make $X$ a left $C$-Galois object, then automatically $C$ will be a Hopf algebra, hence the reflected algebraic quantum group of $A$.\\

\noindent We can again show this by explicitly constructing a homomorphism $\tilde{\beta}_C:C\rightarrow X\otimes X^{\textrm{op}}$, $\beta_C=(\iota\otimes f^{-1})\tilde{\beta}_C$, $\beta_C(c)=c^{[-2]}\otimes c^{[-1]}$. On generators it is given by $\beta_C(u)=x\otimes x^{-1}$ and $\beta_C(w)=y\otimes x^{-m}-1\otimes yx^{-m}$. Again the same chore shows that it has a well-defined extension to $C$, and that it provides the good inverse for the Galois map associated with $\gamma$. This concludes the proof.
\end{proof}

\noindent \emph{Remarks:} 1. If the characteristic of $k$ is zero, then $C$ will not be isomorphic to $A$ when $\mu\neq 0$. For in $A$, the only group-like elements are powers of $a$. Thus any isomorphism would send $u$ to a power $a^l$ of $a$. But then $\mu(a^{lmn}-1)$ would have to be an $n$-th power in $A$, hence, dividing out by $b$, also in $k\lbrack a,a^{-1}\rbrack$. This is impossible.\\
\indent \;\;\;\;\;\;\;\;\;\;\;2. There does not seem to be any straightforward modification of the two-generator Taft algebra Galois objects that provides a Galois object for the \emph{dual} of some $A^{n,m}_\lambda$. It would be interesting to see if such non-trivial Galois objects exist.

\section{Appendix}

\subsection*{\qquad \textit{Multipliers}}

\noindent Let $A$ be a non-degenerate algebra over a field $k$, with or without a unit. The non-degeneracy condition
means that $ab=0$ for all $b\in A$ implies $a=0$, and $ab=0$ for all
$a\in A$ implies $b=0$. As a set, the multiplier algebra $M(A)$ of
$A$ consists of couples $(\lambda,\rho)$, where $\lambda$ and $\rho$
are linear maps $A\rightarrow A$, obeying the following law:
\[a\lambda(b)=\rho(a)b, \qquad\textrm{for all }a,b\in A.\] In
practice, we write $m$ for $(\lambda,\rho)$, and denote $\lambda(a)$
by $ma$ and $\rho(a)$ by $am$. Then the above law is simply an
associativity condition. With the obvious multiplication by
composition of maps, $M(A)$ becomes an algebra, called the
multiplier algebra of $A$. Moreover, if $k=\mathbb{C}$ and $A$ is a $^*$-algebra,
$M(A)$ also carries a $^*$-operation: for $m\in M(A)$ and $a\in A$,
we define $m^*$ by $m^*a=(a^*m)^*$ and $am^*=(ma^*)^*$.\\

\noindent There is a natural map $A\rightarrow M(A)$, letting an element $a$
correspond with left and right multiplication by it. Because of
non-degeneracy, this algebra morphism will be an injection.
In this way, non-degeneracy compensates the possible lack of a unit.
Note that, when $A$ is unital, $M(A)$ is equal to $A$.\\

\noindent Let $B$ be another non-degenerate algebra, and $f$ a non-degenerate
algebra homomorphism $A\rightarrow M(B)$, where by the
non-degeneracy we mean that $f(A)B=B=Bf(A)$. Then
$f$ can be extended to an algebra morphism from $M(A)$ to
$M(B)$, by defining $f(m)(f(a)b)=f(ma)b$ and $(bf(a))f(m)=bf(am)$ for $m\in M(A), a\in A$ and $b\in B$.\\

\subsection*{\qquad \textit{Multiplier Hopf algebras}}

\noindent A \textit{regular multiplier Hopf algebra} consists
of a couple $(A,\Delta)$, with $A$ a non-degenerate algebra,
and $\Delta$, the \textit{comultiplication}, a non-degenerate homomorphism $A\rightarrow M(A\otimes A)$.
Moreover, $(A,\Delta)$ has to satisfy the following conditions:
\begin{itemize}\item[M.1] $(\Delta\otimes \iota)\Delta = (\iota\otimes
\Delta)\Delta$ \qquad (coassociativity). \item[M.2] The maps
\begin{itemize}
\item[]{ $T_{\Delta 2}: A \otimes A \rightarrow M(A \otimes A):
a\otimes b \rightarrow \Delta (a)(1 \otimes b)$, }

\item[]{ $T_{1 \Delta}: A \otimes A \rightarrow M(A \otimes A):
a\otimes b \rightarrow (a \otimes 1) \Delta (b)$, }

\item[]{ $T_{\Delta 1}: A \otimes A \rightarrow M(A \otimes A):
a\otimes b \rightarrow \Delta (a) (b\otimes 1)$, }
\item[]{$T_{2 \Delta}: A \otimes A \rightarrow M(A \otimes A):
a\otimes b \rightarrow (1\otimes a)\Delta (b)$ }
\end{itemize} all induce linear bijections $A\otimes A\rightarrow A\otimes A$.
\end{itemize}

\noindent The $T$-maps can be used to define a counit (which will
be a homomorphism from $A$ to $k$) and an
antipode which will be a linear anti-homomorphism.
Both counit and antipode will be unique, and will satisfy the
corresponding equations of those defining them in the Hopf
algebra case.\\

\noindent When $A$ is a $^*$-algebra over $\mathbb{C}$ and $\Delta$ is $^*$-preserving, we call $(A,\Delta)$ a \textit{regular multiplier Hopf $^*$-algebra}. \\

\subsection*{\qquad \textit{Algebraic quantum groups}}

\noindent An \textit{algebraic quantum group} is a regular multiplier
Hopf algebra $(A,\Delta)$ for which there exists a non-zero
functional $\varphi$ on $A$ such that
\[(\iota\otimes \varphi)(\Delta(a)(b\otimes 1))=\varphi(a)b,\qquad
\textrm{for all } a,b\in A.\] A \textit{$^*$-algebraic quantum group} is an algebraic quantum group which is at the same time a multiplier Hopf $^*$-algebra, such that the functional $\varphi$ is positive:
for every $a\in A$, we have $\varphi(a^*a)\geq 0$. This extra
condition is very restrictive.\\

\noindent For any algebraic quantum group, the functional $\varphi$ will be unique up to multiplication with a
scalar. It will be faithful in the following sense: if
$\varphi(ab)=0$ for all $b\in A$ or $\varphi(ba)=0$ for all $b\in A$, then $a=0$. The functional
$\varphi$ is called the \emph{left invariant functional}. Then
$(A,\Delta)$ will also have a functional $\psi$, such that
\[(\psi\otimes \iota)(\Delta(a)(1\otimes b))=\psi(a)b,
\qquad\textrm{for all } a,b\in A.\] This map $\psi$ is called the
\emph{right invariant functional}. If $A$ is a $^*$-algebraic quantum
group, then $\psi$ can still be chosen so that it is positive.\\

\noindent For any algebraic quantum group, there exists a unique automorphism $\sigma$ of the
algebra $A$, satisfying $\varphi(ab)=\varphi(b\sigma(a))$ for all
$a,b\in A$. It is called \emph{the modular automorphism}. There also exists a unique multiplier $\delta$ such that
\[(\varphi \otimes \iota)(\Delta(a)(1\otimes b))=\varphi(a)\delta
b,\]\[(\varphi\otimes \iota)((1\otimes
b)\Delta(a))=\varphi(a)b\delta,\] for all $a,b\in A$. It is called
\emph{the modular element}.\\

\noindent There is a particular number that can be associated
with an algebraic quantum group. Since $\varphi\circ S^2$ is a left
invariant functional, the uniqueness of $\varphi$ implies there exists
$\tau\in k$ such that $\varphi(S^2(a))=\tau \varphi(a)$, for
all $a\in k$. If $A$ is a $^*$-algebraic quantum group then $\tau=1$. \\

\noindent Note that all formulas which we proved in the second section were known to hold when $X=A$ and $\alpha=\Delta$ (remark that $\beta(a)=S(a_{(1)})\otimes a_{(2)}$ in this case), and we have used some of them in proving our statements.\\

\noindent To any algebraic quantum group $(A,\Delta)$, one can associate another quantum group $(\hat{A},\hat{\Delta})$ which is called its dual. As a set it consists of functionals on $A$ of the form $\varphi(\,\cdot\, a)$ with $a\in A$, where $\varphi$ is the left invariant functional on $A$. Its multiplication and comultiplication are dual to respectively the comultiplication and multiplication on $A$. Intuitively, this means that \[\hat{\Delta}(\omega_1)(a\otimes b)= \omega_1(ab),\]\[(\omega_1\cdot \omega_2)(a)= (\omega_1\otimes \omega_2)(\Delta(a)),\] for $a,b\in A$ and $\omega_1,\omega_2\in \hat{A}$, but some care is needed in giving sense to these formulas.\\

 \noindent The counit on $\hat{A}$ is defined by evaluation in 1, while the antipode is the dual of the antipode of $A$: if $\hat{S}$ denotes the antipode of $\hat{A}$, then \[\hat{S}(\omega_1)(a)=\omega_1(S(a)),\] for $\omega_1\in \hat{A}$ and $a\in A$.
 The left integral $\hat{\varphi}$ of $\hat{A}$ is determined by $\hat{\varphi}(\psi(a\,\cdot\,))=\varepsilon(a)$.\\
  
  \noindent If $A$ is a $^*$-algebraic quantum group, then also $\hat{A}$ will be $^*$-algebraic, with $\hat{\varphi}$ as a positive left invariant functional.

\subsection*{\qquad \textit{Covering issues}}

\noindent When working with multiplier Hopf algebras, it is advantageous to use the Sweedler notation to gain insight into certain formulas. However in this context this is not as straightforward as for Hopf algebras. The problem is that if $(A,\Delta)$ is a multiplier Hopf algebra, then $\Delta(a)$ is an element of $M(A\otimes A)$, and can in general not be written as a sum of elementary tensors. So if we denote $\Delta(a)= \sum a_{(1)}\otimes a_{(2)}$, then this is purely formal, as the right hand side is no well-defined sum of finitely many elements. This gives problems if we want to apply a map to one of the legs of $\Delta(a)$. This is the situation in which we need coverings. For elements of the form $\Delta(a)(1\otimes b)$ with $a,b\in A$ \textit{are} finite sums of elementary tensors in $A\otimes A$, so if we denote this by $\sum a_{(1)}\otimes a_{(2)}b$, there is no trouble in applying a map to the first leg. We then say that `the variable $a_{(2)}$ is covered on the right by $b$'. Although this seems simple, the situation can become quite complicated when multiple coverings are needed (see e.g.\! the examples in \cite{DVW1}). However, in our paper the situation is not so bad, probably because we are working with algebraic quantum groups in stead of general multiplier Hopf algebras: mostly it is seen at first sight if an expression is well-covered or not. This is why we have opted not to emphasize the covering issues too much, since this would probably have obscured certain proofs and statements.








\end{document}